\documentclass[amsfonts]{article}

\usepackage{amssymb,amsmath}
\newtheorem{theorem}{Theorem}[section]

\newtheorem{proposition}[theorem]{Proposition}
\newtheorem{lemma}[theorem]{Lemma}
\newtheorem{definition}[theorem]{Definition}
\newtheorem{remark}[theorem]{Remark}

\def\bR{\mathbb{R}}

\def\bD{\mathbb{D}}

\def\bC{\mathbb{C}}

\def\cA{\mathcal{A}}
\def\cB{\mathcal{B}}

\def\cD{\mathcal{D}}
\def\cE{\mathcal{E}}
\def\cF{\mathcal{F}}
\def\cG{\mathcal{G}}
\def\cH{\mathcal{H}}

\def\cP{\mathcal{P}}

\def\cS{\mathcal{S}}

\def\bxi{\boldsymbol \xi}

\begin{document}

\title{The Stochastic Wave Equation
with Multiplicative Fractional Noise: a Malliavin calculus approach
}

\author{Raluca M. Balan\footnote{Department of Mathematics and Statistics, University of Ottawa,
585 King Edward Avenue, Ottawa, ON, K1N 6N5, Canada. E-mail address:
rbalan@uottawa.ca} \footnote{Research supported by a grant from the
Natural Sciences and Engineering Research Council of Canada.}}

\date{May 28, 2010}

\maketitle

\begin{abstract}
\noindent We consider the stochastic wave equation with
multiplicative noise, which is fractional in time with index
$H>1/2$, and has a homogeneous spatial covariance structure given by
the Riesz kernel of order $\alpha$.
%$f(x)=c_{\alpha,d}|x|^{-(d-\alpha)}$, with
%$0<\alpha<d$.
The solution is interpreted using the Skorohod integral. We show
that the sufficient condition for the existence of the solution is
$\alpha>d-2$, which coincides with the condition obtained in
\cite{dalang99}, when the noise is white in time. Under this
condition, we obtain estimates for the $p$-th moments of the
solution, we deduce its H\"older continuity, and we show that the
solution is Malliavin differentiable of any order. When $d \leq 2$,
we prove that the first-order Malliavin derivative of the solution
satisfies a certain integral equation.
\end{abstract}

{\em MSC 2000 subject classification:} Primary 60H15; secondary
60H07

%60H15 s.p.d.e.
%60H07 stochastic calculus of variations and Malliavin calculus
%60H05 stochastic integrals

{\em Keywords and phrases:} stochastic wave equation,  fractional
Brownian motion, spatially homogeneous Gaussian noise, Malliavin
calculus

\section{Introduction}

In the present article, we consider the following Cauchy problem:
\begin{eqnarray}
\label{wave} \frac{\partial^2 u}{\partial t^2}(t,x)&=& \Delta u
(t,x) +u %\diamond
\dot
W(t,x), \quad t>0, x \in \bR^d \\
\nonumber
u(0, x)&=& u_0(x), \quad x \in \bR^d \\
\nonumber \frac{\partial u}{\partial t}(0,x) &=& u_1(x), \quad x \in
\bR^d,
\end{eqnarray}
where $u_0$ and $u_1$ are deterministic functions, and $W$ is a
zero-mean Gaussian process which is fractional in time, with Hurst
index $H>1/2$, and homogeneous in space, with covariance kernel $f$
(to be specified below). In other words, $W=\{W(\varphi); \varphi
\in \cH \cP\}$ is an isonormal Gaussian process, defined on a
complete probability space $(\Omega,\cF,P)$, with covariance:
$E[W(\varphi)W(\psi)]=\langle \varphi,\psi \rangle_{\cH \cP}$.

Throughout this article, $\cH \cP$ denotes a Hilbert space (which
may contain distributions in $\cS'(\bR^{d+1})$), defined as the
closure of the set $\cE$ of linear combinations of elementary
functions $1_{[0,t] \times A}, t \geq 0, A \in \cB_{b}(\bR^d)$, with
respect to the inner product:
\begin{equation}
\label{def-product-HP} \langle 1_{[0,t] \times A}, 1_{[0,s] \times
B}\rangle_{\cH \cP}:=R_{H}(t,s)\int_{A} \int_{B} f(x-y)dxdy.
\end{equation}
(Here, $\cB_b(\bR^d)$ denotes the class of bounded Borel sets in
$\bR^d$.)

The notation appearing in (\ref{def-product-HP}) needs some
explanation. $R_{H}(t,s)$ denotes the covariance of the fractional
Brownian motion of index $H$, and since we assume that $H>1/2$, we
have:
$$R_{H}(t,s)=\alpha_H \int_0^t \int_0^s |u-v|^{2H-2}dudv,$$
where $\alpha_H=H(2H-1)$. On the other hand, $f:\bR^d \to \bR_{+}$
is the Fourier transform in $\cS'(\bR^d)$, of a tempered measure
$\mu$ on $\bR^d$, i.e.
$$\int_{\bR^d} f(x)\phi(x)dx=\int_{\bR^d} \cF \phi(\xi)\mu(d\xi), \quad \forall \phi \in \cS(\bR^d),$$
where $\cF \phi(\xi)=\int_{\bR^d}e^{-i \xi \cdot x}\phi(x)dx$ is the
Fourier transform of $\phi$, and $\cS(\bR^d)$ is the space of
rapidly decreasing $C^{\infty}$-functions on $\bR^d$.

In the present article, a solution of (\ref{wave}) is an adapted
square-integrable process $u=\{u(t,x);t \geq 0, x \in \bR^d\}$ which
satisfies the following integral equation:
\begin{equation}
\label{u-integral-eq} u(t,x)=w(t,x)+\int_0^t \int_{\bR^d}
G(t-s,x-y)u(s,y)W(\delta s, \delta y),
\end{equation}
where $G$ is the fundamental solution of the wave operator,
%$u_{tt}-\Delta u=0$,
the stochastic integral is interpreted in the Skorohod sense, %(i.e. it is defined using the divergence operator $\delta$, as explained in Section %\ref{framework} below),
and $w$ is the solution of the equation $w_{tt}=\Delta w$, with
initial conditions $w(0,\cdot)=u_0$, $w_t(0,\cdot)=u_1$.
%\begin{eqnarray*}
%\label{eq-w}
%\frac{\partial^2 w}{\partial t^2}(t,x)&=& \Delta w (t,x) \quad t>0, x \in \bR^d \\
%\nonumber
%w(0, x)&=& u_0, \quad x \in \bR^d \\
%\nonumber \frac{\partial w}{\partial t}(0,x) &=& u_1, \quad x \in \bR^d.
%\end{eqnarray*}
%(see p.11 of \cite{DMT08}, for the case $H=1/2$).

We are interested in the case $u_0=1$ and $u_1=0$, and hence
$w(t,x)=1$ for all $t \geq 0,x \in \bR^d$. (In the case $H=1/2$,
this corresponds to the equation $u_{tt}=\Delta u+(u+1)\dot W$ with
zero initial conditions.)

The case of the heat equation with multiplicative fractional noise
was treated in \cite{hu01}, \cite{hu-nualart09} and
\cite{balan-tudor09},
%(in the case $f(x)=\prod_{i=1}^{n}\alpha_{H_i}|x_i|^{2H_i-2}$ with $1/2<H_i<1$, $f(x)=\delta_0$, %respectively $f(x)=c_{\alpha,d}|x|^{-(d-\alpha)}$ with $0<\alpha<d$) ,
while the recent article \cite{hu-nualart-song09} gives a
Feynman-Kac representation for the solution. In the case of the heat
equation,
%when $f$ is the Riesz kernel $f(x)=c_{\alpha,d}|x|^{-(d-\alpha)}$, with $0<\alpha<d$,
the key estimate which leads to a sufficient condition for the
existence of the solution is:
\begin{equation}
\label{elementary-est-heat-intro}\int_{\bR^d}
e^{-t|\xi|^2/2}|\xi-\eta|^{-\alpha}d\xi \leq C_{\alpha,d}
t^{-(d-\alpha)/2}, \quad \mbox{for any} \ t>0,\eta \in \bR^d,
\end{equation}
for any $0<\alpha<d$ (see Lemma 6.1 of \cite{hu01}, or Lemma 3.3 of
\cite{balan-tudor09}).
%A similar estimate is given by Lemma 6.1 of \cite{hu01}, in the case %$f(x)=\prod_{i=1}^{n}\alpha_{H_i}|x_i|^{2H_i-2}$ with $1/2<H_i<1$.

The following estimate lies at the origin of our developments, being
the analogue of (\ref{elementary-est-heat-intro}), which is needed
for the wave equation:
$$\int_{\bR^d}\frac{\sin^2(t|\xi|)}{|\xi|^2}|\xi-\eta|^{-\alpha}d\xi\leq
C_{\alpha,d}t^{\alpha-d+2}, \quad \mbox{for any} \ t>0,\eta \in
\bR^d,$$ for any $d-2<\alpha<d$ (see Lemma \ref{hu-lemma} below).

As far as we know, the only other study of the wave equation, in which the fractional noise enters the equation in a multiplicative way is \cite{quer-tindel06}, which treats the case $d=1$. %and $\sigma(u)$ instead of $u$).
The authors of \cite{quer-tindel06} use a pathwise integral for
interpreting the solution, instead of the Skorohod integral used in
the present article, which makes it difficult to compare the
results.

The study of the wave and heat equations driven by a Gaussian noise
which is white in time and homogeneous in space was considered by
many authors, using martingale methods (see \cite{dalang99},
\cite{millet-sanzsole99}, \cite{sanzsole05}, \cite{conus-dalang09},
and the references therein). These methods work for more general
equations (in which the factor $u$ multiplying the noise may be
replaced by $\sigma(u)$, for a Lipschitz function $\sigma$), but
fail in the fractional case. The method that we use in the present
article is specific to the case $\sigma(u)=u$, in which the solution
has a Wiener chaos decomposition whose kernels can be written down
in closed form.

The fact that the solution has a known Wiener chaos decomposition
has several implications. More precisely, once we have a uniform
bound for the moments of order $2$ of the solution, we can pass to
the moments of order $p > 2$, using the hypercontractivity of the
Ornstein-Uhlenbeck semigroup. In this manner, we obtain %almost immediately
estimates for the $p$-th moments of the increments of the solution
(which yield its H\"older continuity),
%(which yield the existence of a H\"older continuous modification),
and we show %with little effort,
the Malliavin differentiability of the solution. % (of any order).
Finally, assuming that $d \leq 2$, we prove that the first-order
Malliavin derivative of the solution satisfies an integral equation,
using a
Hilbert-space-valued Skorohod integral. %, which seems to be a new idea.
These results are valid for the heat equation as well.
%This article is organized as follows. In Section \ref{framework}, ... In Section \ref{existence},
%we identify a sufficient condition for the existence of the solution in the case of particular kernels $f$ %(e.g. Riesz kernel). In Section ...

The present article leaves some open problems:

1. For the existence of the solution, we could treat only the case
of the Riesz kernel $f(x)=c_{\alpha,d}|x|^{-(d-\alpha)}$ with
$0<\alpha<d$, but other particular cases (e.g.
$f(x)=\prod_{i=1}^{d}\alpha_{H_i}|x_i|^{2H_i-2}$ with $1/2<H_i<1$)
would be similar. However, we do not know what is the sufficient
condition for the existence of the solution, in the case of a {\em
general} kernel $f$. Even in the case of the Riesz kernel, it is not
clear if the
sufficient condition $\alpha>d-2$ %(which does not depend on $H$)
is {\em necessary} for the existence of the solution. These
questions are open in the case of the heat equation as well.
%The necessity question is open for $H=1/2$ too.

%We do not know if the condition $\alpha>d-2$ is necessary for the existence
%of the solution. It is weird that this condition does not depend on $H$
%(indeed this is the condition encountered in the case of the white noise in
%time). We believe that the necessary and sufficient condition is stronger
%than $\alpha>d-2H-1$ (encountered in the case case of additive noise,
%studied in \cite{BT10}), but maybe weaker %than $\alpha>d-2$.

3. In the case of the white noise in time (see e.g.
\cite{nualart-quer07}), the argument which leads to the existence of
the density of $u(t,x)$ relies on the fact that $Du(t,x)$ satisfies
a certain integral equation, and its $\cH \cP$-norm coincides with
the norm in $L^2((0,t); \cH)$, for a certain Hilbert space $\cH$ of
distributions in $\cS'(\bR^d)$. Therefore, in the calculation of
$\|Du(t,x)\|_{\cH \cP}$, one can bound from below the integral over
$(0,t)$ with the integral over a small interval $(t,t-\delta)$. This
bound is not justified in the fractional case.
 As we could not find an alternative argument to prove that
 $\|Du(t,x)\|_{\cH \cP}>0$ a.s.
 (which would allow us to apply the Hirsch-Bouleau criterion),
 we could not show that the law of $u(t,x)$ has a density, even if
 $d \leq 2$. %is absolutely continuous with respect to the Lebesgue measure.% on $\bR$.

\section{Preliminary Results}
\label{framework}

%The material presented in this section is usually taken for granted by
% the specialists, who can skip reading it and move directly to
% Section \ref{existence}. However, we feel that the article
%would not be complete without this discussion. %We include it for the sake of completeness.

Intuitively, %as in the ,% the kernels $f_n(\cdot, t,x)$ appearing in the Wiener chaos decomposition of %$u(t,x)$ should be written in closed form using $G$. More precisely,
the solution of (\ref{wave}) should be given by a series of iterated
integrals:
$$u(t,x)=1+\int_0^t \int_{\bR^d} G(t-t_1,x-x_1) W(dt_1,dx_1)+ $$
$$\int_0^t \int_{\bR^d} G(t-t_2,x-x_2)\left(\int_0^{t_2} \int_{\bR^2}G(t_2-t_1,x_2-x_1)W(dt_1,dx_1)\right)W(dt_2,dx_2)+\ldots$$

Since in dimension $d \geq 3$, $G(t,\cdot)$ is a distribution in
$\bR^d$, the product
$$f_n(t_1,x_1, \ldots, t_n,x_n,t,x):=G(t-t_n,x-x_n)G(t_n-t_{n-1},x_n-x_{n-1}) \ldots $$
\begin{equation}
\label{formal-def-fn}G(t_2-t_1,x_2-x_1)1_{\{ 0<t_1< \ldots<t_n<t\}},
\end{equation}
has to be defined as the product of distributions, and one has to be
careful with the definition of the iterated integrals above.

The goal of this section is to take care of this difficulty, by
tackling the following three problems:
\begin{itemize}
\item (Subsection \ref{def-kernels}) For any $t>0,x \in \bR^d$ and for any $0<t_1< \ldots<t_n<t$, we give a meaning to $f_n(t_1, \cdot, \ldots, t_n, \cdot,t,x)$ as a distribution in $\cS'(\bR^{nd})$, and calculate its Fourier transform.

 \item (Subsection \ref{def-space-HP-n}) For any $t>0, x \in \bR^d$, we give a general criterion which ensures that $f_n(\cdot,t,x) \in \cH \cP^{\otimes n}$.
     %if $\|f_n(\cdot, t,x)\|_{\cH \cP^{\otimes n}}^2<\infty$.

 \item (Subsection \ref{Malliavin}) %Let $u(0,x)=1$ for any $x \in \bR^d$.
 Suppose that for any $t>0,x \in \bR^d$, $f_n(\cdot, t,x) \in
 \cH \cP^{\otimes n}$, and the series
 $u(t,x):=1+\sum_{n \geq 1}I_n(f_n(\cdot,t,x))$
converges in $L^2(\Omega)$, where $I_n$ denotes the multiple Wiener
integral with respect to $W$. We show that $u=\{u(t,x); t \geq 0, x
\in \bR^d\}$ is a %(random-field)
solution of (\ref{wave}), in a sense which will be described below.
\end{itemize}

\subsection{The definition of the kernels $f_n(\cdot,t,x)$}
\label{def-kernels}

%We begin with item 1 above.

In this subsection, we give a rigorous meaning to the kernels
$f_n(\cdot,t,x)$. We let $C_{0}^{\infty}(\bR^d)$ be the space of
infinitely differentiable functions on $\bR^d$ with compact support,
and $\cD'(\bR^d)$ be the space of (Schwartz) distributions on
$\bR^d$.

Assume first that $n=2$ and let $0<t_1<t_2<t$ be arbitrary. We
proceed to the formal calculation of the action of
$f_2(t_1,\cdot,t_2, \cdot, t,x)$ on a test function $\phi=\phi_1
\otimes \phi_2$ with $\phi_1, \phi_2 \in C_0^{\infty}(\bR^d)$:
\begin{eqnarray*}
(f_2(t_1,\cdot,t_2,\cdot,t,x), \phi)&=&\int_{\bR^d} G(t-t_2,x-x_2)\phi_2(x_2) \int_{\bR^d} G(t_2-t_1,x_2-x_1)\phi_1(x_1)dx_1 dx_2 \\
&=& \int_{\bR^d} G(t-t_2,x-x_2)\phi_2(x_2) \varphi_1(t_2-t_1,x_2) dx_2 \\
&=&  \int_{\bR^d} G(t-t_2,x-x_2)\psi_2(t_2-t_1,x_2) dx_2 \\
&=&
%(\psi_2(t_2-t_1,\cdot)*G(t-t_2,\cdot))(x)=
\varphi_2(t_2-t_1,t-t_2,x)
\end{eqnarray*}
where
\begin{eqnarray*}
& & \psi_1(\cdot)= \phi_1(\cdot), \quad \quad \quad \quad \quad \mbox{and} \quad \psi_2(s,\cdot)=\phi_2(\cdot)\varphi_1(s,\cdot) \\
& & \varphi_1(s,\cdot)= \psi_1(\cdot) * G(s,\cdot) , \quad \quad
\quad \varphi_2(s_1,s_2,\cdot)=\psi_2(s_1,\cdot)*G(s_2,\cdot),
\end{eqnarray*}
and $*$ denotes the convolution with respect to the space variable. Similar formal calculations can be done for any $n$. %$(f_n(t_1, \ldots, \cdots,t_n,x_n,t,x),\phi)$.

Based on these calculations, for any $0<t_1<\ldots<t_n<t$ fixed, we
let $f_n(t_1,\cdot,\ldots, t_n,\cdot, t,x)$ be the element of
$\cD'(\bR^{nd})$ whose action on a test function $\phi=\phi_1
\otimes \ldots \otimes \phi_n$ with $\phi_i \in
C_0^{\infty}(\bR^d)$, is given by:
\begin{equation}
\label{definition-fn-in-S'} (f_n(t_1,\cdot,\ldots, t_n,\cdot, t,x),
\phi):= \varphi_n(t_2-t_1,t_3-t_2,\ldots,t-t_n,x),
\end{equation}
where the pairs $(\psi_k,\varphi_k)$ are defined recursively for
$k=1,\ldots,n$ by the following relations:
\begin{eqnarray}
\label{def-psi}
\psi_k(s_1,\ldots,s_{k-1},\cdot)&=& \phi_k(\cdot)\varphi_{k-1}(s_1, \ldots,s_{k-1},\cdot) \\
\label{def-varphi} \varphi_k(s_1,\ldots,s_k,\cdot)&=&
\psi_k(s_1,\ldots,s_{k-1},\cdot)*G(s_k,\cdot).
\end{eqnarray}
Note that $\psi_k(s_1, \ldots, s_{k-1},\cdot) \in
C_0^{\infty}(\bR^d)$ and $\varphi_k(s_1, \ldots,s_k,\cdot) \in
\cS(\bR^d)$, since $G(s, \cdot)$ is a distribution with rapid
decrease in $\bR^d$ (see p. 245 of \cite{schwartz66}). The previous
definition is extended to $\phi=\phi_1 \otimes \ldots \otimes
\phi_n$ with $\phi_i \in \cS(\bR^d)$.

%{\bf We have $f_n(t_1, \cdot, \ldots, t_n,\cdot,t,x) \in \cS'(\bR^{nd})$. }

The next result shows that the Fourier transform of $f_n(t_1, \cdot,
\ldots, t_n, \cdot, t,x)$ is a function in $\bR^{nd}$, given by:
\begin{eqnarray}
\nonumber \cF f_n(t_1, \cdot, \ldots, t_n, \cdot, t,x)(\xi_1,
\ldots, \xi_n)=
e^{-i(\xi_1+\ldots+\xi_n)\cdot x}\overline{\cF G(t_2-t_1, \cdot)(\xi_1)} \\
\label{Fourier-fn-formula} \overline{\cF G(t_3-t_2,
\cdot)(\xi_1+\xi_2)}\ldots  \overline{\cF G(t-t_n,
\cdot)(\xi_1+\ldots+\xi_n)}.
\end{eqnarray}

%The proof of Proposition \ref{Fourier-fn-prop} can be found in appendix.

\begin{proposition}
\label{Fourier-fn-prop} For any $0<t_1< \ldots<t_n<t$ and for any
$h=h_1 \otimes \ldots \otimes h_n$ with $h_i \in
C_0^{\infty}(\bR^d)$, we have:
\begin{eqnarray*}
\lefteqn{(f_n(t_1,\cdot,\ldots, t_n,\cdot, t,x), h)=\int_{\bR^{nd}}h(\xi_1, \ldots,\xi_n)e^{-i(\xi_1+\ldots+\xi_n)\cdot x} \overline{\cF G(t_2-t_1, \cdot)(\xi_1)} }\\
 & & \overline{\cF G(t_3-t_2, \cdot)(\xi_1+\xi_2)} \ldots  \overline{\cF G(t-t_n, \cdot)(\xi_1+\ldots+\xi_n)}d\xi_1 \ldots d\xi_n.
\end{eqnarray*}
\end{proposition}

\noindent {\bf Proof:} Note that $\phi:=\cF h=\phi_1 \otimes \ldots
\otimes \phi_n$, where $\phi_i:=\cF h_i \in \cS(\bR^d)$. By the
definition of the Fourier transform in $\cS'(\bR^d)$ and
(\ref{definition-fn-in-S'}), we have:
\begin{eqnarray}
\nonumber
(\cF f_n(t_1, \cdot, \ldots, t_n, \cdot, t,x),h)&=&(f_n(t_1, \cdot, \ldots, t_n, \cdot, t,x),\phi) \\
\label{Fourier-fn}&=& \varphi_n(t_2-t_1,t_3-t_2,\ldots,t-t_n,x),
\end{eqnarray}
where $(\psi_k,\varphi_k), 1 \leq k \leq n$ are defined recursively
by (\ref{def-psi})-(\ref{def-varphi}).

We proceed to the evaluation of
$\varphi_n(t_2-t_1,t_3-t_2,\ldots,t-t_n,x)$.

{\em Step 1.} For any $s_1, \ldots, s_{n-1} \in [0,t]$, we define
recursively the following functions: $g_1=h_1$,
$$g_k(s_1, \ldots, s_{k-1}, \cdot)=h_k * (g_{k-1}(s_1, \ldots, s_{k-2},\cdot) \overline{\cF G(s_{k-1}, \cdot)}), \quad 2 \leq k \leq n.$$

By induction on $k$, $2 \leq k \leq n$, it follows that $g_k(s_1,
\ldots, s_{k-1}, \cdot) \in C_0^{\infty}(\bR^d)$ (since $\cF G (s,
\cdot)$ is a $C^{\infty}$-function on $\bR^d$), and
\begin{equation}
\label{new-def-gk} g_{k}(s_1, \ldots, s_{k-1},
\eta_k)=\int_{\bR^{(k-1)d}} h_1(\eta_1) h_2(\eta_2-\eta_1) \ldots
h_k(\eta_k-\eta_{k-1}) \overline{\cF G(s_1, \cdot)(\eta_1)} \ldots
$$
$$\overline{\cF G(s_{k-1}, \cdot)(\eta_{k-1})} d\eta_{1} \ldots d\eta_{k-1}.
\end{equation}

%For $k=2$, we have: $$g_2(s,\eta_2)=[h_2* (h_1\overline{\cF G(s,\cdot)})](\eta_2)=\int_{\bR^d} %h_2(\eta_2-\eta_1)h_1(\eta_1) \overline{\cF G(s,\cdot)(\eta_1)}d\eta_1.$$

{\em Step 2.} We prove by induction on $k$, $1 \leq k \leq n$, that:
\begin{equation}
\label{new-def-varphi} \varphi_k(s_1, \ldots, s_k, \cdot)=\cF
[g_{k}(s_1, \ldots, s_{k-1}, \cdot) \overline{\cF G(s_k, \cdot)}].
\end{equation}

Before proving (\ref{new-def-varphi}), note that:
\begin{equation}
\label{Fourier-identity} \cF g * G(s, \cdot)=\cF (g \overline{\cF
G(s, \cdot)}), \quad \forall  g \in C_0^{\infty}(\bR^d), \forall
s>0,
\end{equation}
since $(\cF g) * G(s,\cdot)=\cF [\cF^{-1}(\cF g * G(s,\cdot))]=
%\cF [\cF^{-1}(\cF g) \cF^{-1}G(s,\cdot)]=
\cF(g \overline{\cF G(s, \cdot)})$.

For $k=1$, we have $g_1=h_1=\cF^{-1} \phi_1$ and
$$\varphi_1(s, \cdot)=\phi_1 * G(s,\cdot)=\cF g_1 * G(s,\cdot)=\cF (g_1 \overline{\cF G(s, \cdot)}),$$
where we used (\ref{Fourier-identity}) for the last equality. This
proves (\ref{new-def-varphi}) for $k=1$.

Suppose that (\ref{new-def-varphi}) holds for $k-1$. Then
\begin{eqnarray*}
\psi_k(s_1, \ldots, s_{k-1},x)&=&\phi_k(x) \varphi_{k-1}(s_1, \ldots, s_{k-1}, x) \\
&=& \cF h_k(x) \cF[g_{k-1}(s_1, \ldots, s_{k-2}, \cdot) \overline{\cF G(s_{k-1},\cdot)}](x)\\
&=& \cF [h_k *(g_{k-1}(s_1, \ldots, s_{k-2}, \cdot) \overline{\cF G(s_{k-1},\cdot)})](x)\\
&=& \cF g_{k}(s_1, \ldots, s_{k-1}, \cdot)(x)
\end{eqnarray*}
and
\begin{eqnarray*}
\varphi_k(s_1, \ldots, s_k,\cdot)&=& \psi_k(s_1, \ldots, s_{k-1},
\cdot)*G(s_{k}, \cdot)=
\cF g_{k}(s_1, \ldots, s_{k-1}, \cdot)*G(s_k,\cdot)\\
&=& \cF[g_k(s_1, \ldots,s_{k-1}, \cdot) \overline{\cF G(s_k,
\cdot)}],
\end{eqnarray*}
where we used (\ref{Fourier-identity}) for the last equality. This
proves (\ref{new-def-varphi}).

{\em Step 3.} Using (\ref{new-def-varphi}) and (\ref{new-def-gk}),
we obtain that for any $1 \leq k \leq n$,
\begin{eqnarray}
\nonumber
\lefteqn{\varphi_k(s_1, \ldots,s_{k},x)=\int_{\bR^d} e^{-i \eta_k \cdot x} g_k(s_1, \ldots,s_{k-1},\eta_k)\overline{\cF G(s_k, \cdot)(\eta_k)}d\eta_k } \\
\nonumber && = \int_{\bR^{kd}} e^{-i \eta_k \cdot x} h_1(\eta_1) h_2(\eta_2-\eta_1) \ldots h_k(\eta_k-\eta_{k-1}) \overline{\cF G(s_1, \cdot)(\eta_1)} \\
\nonumber & & \quad \ldots \overline{\cF G(s_{k-1}, \cdot)(\eta_{k-1})}\ \overline{\cF G(s_k, \cdot)(\eta_k)} d\eta_{1} \ldots d\eta_k \\
\label{def-varphi-k} & & = \int_{\bR^{kd}}e^{-i(\xi_1+\ldots+\xi_k)
\cdot x} \prod_{j=1}^{k}h_j(\xi_j) \prod_{j=1}^{k}\overline{\cF
G(s_j, \cdot)(\xi_1+\ldots+ \xi_j)}  d\xi_1 \ldots d\xi_k,
\end{eqnarray}
where for the last equality we used the change of variables
$\xi_1=\eta_1$, $\xi_j=\eta_j-\eta_{j-1}$, for $2 \leq j \leq k$. We
now use (\ref{Fourier-fn}). The conclusion follows using
(\ref{def-varphi-k}) with $k=n$,  $s_n=t-t_{n}$ and
$s_i=t_{i+1}-t_i$ for $1 \leq i \leq n-1$. $\Box$

\subsection{The space $\cH \cP^{\otimes n}$}
\label{def-space-HP-n}

%We continue with item 2 above. For this,

In this subsection, we give a criterion for checking that an element
$\varphi \in \cH \cP^{\otimes n}$, which can be viewed as a
multi-dimensional analogue of Theorem 2.1 of \cite{BT10}. This
criterion is then applied to the case $\varphi=f_n(\cdot,t,x)$.

For this purpose, for any $T>0$ fixed, we define the
multi-dimensional transfer operator $K_{H,n}^*$ by:
$$(K_{H,n}^* 1_{[0,t_1] \times \ldots [0,t_n]})(s_1, \ldots,s_n):=\prod_{i=1}^{n}(K_{H}^* 1_{[0,t_i]})(s_i), \quad s_1, \ldots,s_n \in (0,T).$$

Let $\cE_{\bC}(0,T)$ be the set of all complex linear combinations
of indicator functions $1_{[0,t]},t \in [0,T]$, and $\cH_{\bC}(0,T)$
be the closure of $\cE_{\bC}(0,T)$  with respect to the inner
product:
$$\langle \varphi, \psi \rangle_{\cH_{\bC}(0,T)}=\alpha_H \int_0^T
\int_O^T \varphi(u) \overline{\psi(v)}|u-v|^{2H-2}du dv.$$

The operator $K_{H,n}^*$ is an isometry between
$\cE_{\bC}(0,T)^{\otimes n}$ and $L^2((0,T)^n)$, which can be
extended to $\cH_{\bC}(0,T)^{\otimes n}$. In terms of fractional
integrals, we have:
$$(K_{H,n}^* \phi)({\bf s})=(c_{H}^*)^n \Gamma(H-1/2)^n [{\bf s}]^{1/2-H} I_{T-,n}^{H-1/2}([{\bf u}]^{H-1/2} \phi({\bf u}))({\bf s}),$$
where $c_{H}^*=(\frac{\alpha_H}{\beta(H-1/2,2-2H)})^{1/2}$, ${\bf
s}=(s_1, \ldots,s_n)$, $[{\bf s}]= s_1 \ldots s_n$, and
$$I_{T-,n}^{\alpha} f ({\bf s}):=\frac{1}{\Gamma(\alpha)^n} \int_{s_1}^{T} \ldots \int_{s_n}^{T} [{\bf u}-{\bf s}]^{\alpha-1}f({\bf u})d{\bf u}$$
is a multi-dimensional fractional integral of $f \in L^1((0,T)^n)$,
of order $\alpha \in (0,1)$.

For any function $\phi \in \cH_{\bC}(0,T)^{\otimes n}$, we have:
\begin{eqnarray}
\nonumber \lefteqn{\alpha_{H}^{n} \int_{(0,T)^{2n}}\phi({\bf u}) \overline{\phi({\bf v})}[{\bf u}-{\bf v}]^{2H-2}d{\bf u}d{\bf v}=}\\
\label{proof-phi-in-HP-1} & & d_{H}^n
\int_{(0,T)^n}|I_{T-,n}^{H-1/2}([{\bf u}]^{H-1/2} \phi({\bf
u}))({\bf s})|^2 \lambda_{H,n}(d{\bf s}),
\end{eqnarray}
where $d_{H}=(c_{H}^*)^2 \Gamma(H-1/2)^2$ and $\lambda_{H,n}(d{\bf
s})=[{\bf s}]^{1-2H}d{\bf s}$.

Let $\cE_T$ be the class of elementary functions on $(0,T) \times
\bR^d$. For any $\varphi \in \cE_T^{\otimes n}$,
\begin{eqnarray}
\nonumber
\|\varphi\|_{\cH \cP^{\otimes n}}^{2}&=&d_{H}^{n}\int_{\bR^{nd}} \int_{(0,T)^n}|I_{T-,n}^{H-1/2}([{\bf u}]^{H-1/2} \cF \varphi(u_1, \cdot, \ldots,u_n,\cdot)(\bxi))({\bf s})|^2 \\
\label{proof-phi-in-HP-2} & & \lambda_{H,n}(d{\bf s}) \mu(d\xi_1)
\ldots \mu(d\xi_n),
\end{eqnarray}
where $\bxi=(\xi_1, \ldots,\xi_n) \in \bR^{nd}$.

The following theorem is the multidimensional analogue of Theorem
2.1 of \cite{BT10}.

\begin{theorem}
\label{th-varphi-in-HP}Let $(0,T)^n \ni {\bf t} \mapsto \varphi(t_1,
\cdot, \ldots, t_n,\cdot) \in \cS'(\bR^{nd})$ be a deterministic
function such that $\cF f(t_1, \cdot, \ldots, t_n \cdot)$ is a
function on $\bR^{nd}$, for all ${\bf t} \in (0,T)^n$.

Suppose that:\\
(i) the function ${\bf t} \mapsto \cF \varphi(t_1, \cdot, \ldots,t_n, \cdot)(\bxi)$ belongs to $\cH_{\bC}(0,T)^{\otimes n}$ for all $\bxi\in \bR^{nd}$;\\
(ii) the function $({\bf t},\bxi) \mapsto \cF \varphi(t_1, \cdot, \ldots,t_n, \cdot)(\bxi)$ is measurable on $(0,T)^n \times \bR^{nd}$;\\
(iii) $\int_{(0,T)^n} \prod_{i=1}^{n}1_{\{u_i \geq s_i\}}[{\bf
u}]^{H-1/2}[{\bf u}-{\bf s}]^{H-3/2}
|\cF \varphi(u_1, \cdot, \ldots,u_n, \cdot)(\bxi)|d{\bf u}<\infty$ for all $({\bf s},\bxi) \in (0,T)^n \times \bR^{nd}$ (or $\cF \varphi(s_1, \cdot, \ldots,s_n, \cdot)(\bxi) \geq 0$ for all $({\bf s}, \bxi)$). %\in (0,T)^n \times \bR^{nd}$).

If
$$I_T:=\alpha_H^n \int_{\bR^{nd}} \int_{(0,T)^{2n}} \cF \varphi(u_1, \cdot, \ldots,u_n, \cdot)(\bxi)\overline{\cF \varphi(v_1, \cdot, \ldots,v_n, \cdot)(\bxi)}$$
\begin{equation}
\label{IT-finite} \prod_{i=1}^{n}|u_i-v_i|^{2H-2} d{\bf u}d{\bf v}
\mu(d\xi_1) \ldots \mu(d \xi_n)<\infty,
\end{equation}
then $\varphi \in \cH \cP^{\otimes n}$ and $\|\varphi\|_{\cH
\cP^{\otimes n}}^{2}=I_T$. (By convention, we set $\varphi(t_1,
\cdot, \ldots, t_n, \cdot)=0$ if $t_i>T$ for some $i=1, \ldots,n$.)
\end{theorem}

\noindent {\bf Proof:}  The argument is similar to the one used in
the proof of Theorem 2.1 of \cite{BT10}, being based on relations
(\ref{proof-phi-in-HP-1}) and (\ref{proof-phi-in-HP-2}) above. We
omit the details. $\Box$
%(applied to the function $\phi_{\bxi}({\bf t})=\cF \varphi(t_1, \cdot, \ldots, t_n,\cdot)(\bxi)$, for any $\bxi \in \bR^{nd}$ fixed),

\vspace{3mm}

\begin{remark}
\label{rem-varphi-in-HP} {\rm In our case, we apply Theorem
\ref{th-varphi-in-HP} to the function $$(0,t)^n \ni {\bf t} \mapsto
\varphi(t_1, \cdot, \ldots, t_n,\cdot)=f_n(t_1, \cdot, \ldots,
t_n,\cdot,t,x).$$ (We define $f_n(t_1, \cdot, \ldots,
t_n,\cdot,t,x)$ to be $0$ if the relation $t_1<\ldots<t_n$ is not
satisfied.) To see that this function satisfies hypothesis (i)-(iii)
of Theorem \ref{th-varphi-in-HP}, we use (\ref{Fourier-fn-formula})
and the fact that:
\begin{equation}
\label{FG-bounded} |\cF G(t,\cdot)(\xi)| \leq C_T \quad \forall t
\in [0,T], \xi \in \bR^d,
\end{equation}
From here, we infer that the map ${\bf t} \mapsto \cF f_n(t_1,
\cdot, \ldots,t_n,\cdot,t,x)$ belongs to $L^2_{\bC}(0,T)^{\otimes
n}$, which is included in $\cH_{\bC}(0,T)^{\otimes n}$.

Therefore, to show that $f_n(\cdot,t,x) \in \cH \cP^{\otimes n}$, it
suffices to prove that (\ref{IT-finite}) holds. This will be done in
Section \ref{existence}.}
\end{remark}

\subsection{Malliavin Calculus}
\label{Malliavin}

%We now tackle item 3 mentioned at the beginning of this section. For this,

In this subsection, we introduce the basic elements of the Malliavin
calculus with respect to the isonormal Gaussian process $W$ (see
\cite{nualart06} for more details).

We first introduce the multiple Wiener integral with respect to $W$.
%linear subspace of $L^2(\Omega)$ generated by the random variables $\{H_n(W(\varphi)); \varphi \in \cH %\cP, \|\varphi \|_{\cH \cP}=1\}$, where $H_n(x)$ is the $n$-th Hermite polynomial. Let $\cH \cP_0=\bR$.
Let $\cG$ be the $\sigma$-field generated by $\{W(\varphi); \varphi
\in \cH \cP\}$. By Theorem 1.1.1 of \cite{nualart06}, $L^2(\Omega,
\cG,P)=\oplus_{n=0}^{\infty}\cH \cP_n$, where $\cH \cP_n$ be the
$n$-th Wiener chaos of $W$. Hence, every $F \in L^{2}(\Omega,\cG,P)$
admits the following Wiener chaos expansion:
\begin{equation}
\label{Wiener-chaos-0} F=\sum_{n=0}^{\infty}J_nF,
\end{equation} where $J_n$ is the projection on $\cH \cP_n$, for $n \geq 1$. By convention, $J_0 F=E(F)$.

The definition of the multiple Wiener integral $I_n$ is similar to the white noise case (see Subsection 1.1.2 of \cite{nualart06}). %We include this definition for the sake of completeness.
More precisely, $I_n$ is a linear and continuous operator from $\cH
\cP^{\otimes n}$ onto $\cH \cP_n$. For any $f \in \cH \cP^{\otimes
n}$, we write
$$I_n(f)=\int_{(\bR_{+} \times \bR^{d})^n}f(t_1,x_1,\ldots,
t_n,x_n)W(dt_1,dx_1) \ldots W(dt_n, dx_n),$$ even if $f_n$ is not a
function in $(t_1,x_1), \ldots , (t_n,x_n)$. Note that
$I_n(f)=I_n(\tilde f)$ and
\begin{equation}
\label{isometry-In} E(I_n(f) I_{n}(g))=n! \langle \tilde f, \tilde g
\rangle_{\cH \cP^{\otimes n}},
 \end{equation}
 where $\tilde f$ denotes the symmetrization of $f$, i.e.
$$\tilde f (t_1,x_1, \ldots,t_n,x_n)=\frac{1}{n!}\sum_{\rho \in S_n}f(t_{\rho(1)},x_{\rho(1)}, \ldots,
t_{\rho(n)},x_{\rho(n)}).$$ (Here $S_n$ denotes the set of all
permutations of $\{1, \ldots,n\}$.)
%For any $\varphi \in \cH \cP$ with $\|\varphi \|_{\cH
%\cP}=1$,
%$$I_n(\varphi^{\otimes n})=n! \ H_n(W(\varphi)),$$ and
%$I_n$ maps $\cH \cP^{\otimes n}$ onto $\cH \cP_n$. Hence,

Any random variable $F \in L^2(\Omega,\cG,P)$ admits the
decomposition:
\begin{equation}
\label{Wiener-chaos-F}
F=\sum_{n \geq 0} I_n(f_n)%=\sum_{n \geq 0} I_n(\tilde f_n),
\end{equation}
for some $f_n \in \cH \cP^{\otimes n}$ symmetric, and  %Moreover, we may assume that $(f_n)_n$ are symmetric (i.e. $f_n =\tilde f_n$), %and in this case, they are uniquely determined by $F$.
$$E|F|^2=\sum_{n=0}^{\infty}
E|I_n(f_n)|^2=\sum_{n=0}^{\infty}n! \ \|f_n\|_{\cH \cP^{\otimes
n}}^2.$$ (By convention, $f_0=E(F)$ and $I_0(x)=x$ for all $x \in
\bR$.)

We now introduce the derivative operator. Let $\cS$ be the class of
smooth random variables of the form
\begin{equation}
\label{form-F}F=f(W(\varphi_1),\ldots, W(\varphi_n)),
\end{equation} where $f \in C_{b}^{\infty}(\bR^n)$, $\varphi_i \in \cH \cP$, $n \geq 1$, and
$C_b^{\infty}(\bR^n)$ is the class of bounded $C^{\infty}$-functions
on $\bR^n$, whose partial derivatives are bounded. The Malliavin
derivative of $F$ of the form (\ref{form-F}) is an $\cH \cP$-valued
random variable given by:
$$DF:=\sum_{i=1}^{n}\frac{\partial f}{\partial x_i}(W(\varphi_1),\ldots,
W(\varphi_n))\varphi_i.$$
%Note that $DF \in L^2(\Omega; \cH \cP)$. By abuse of
%notation, we sometimes write $DF=\{D_{t,x} F; t \geq 0, x \in \bR^d\}$
%even if $D_{t,x} F$ is not a function in $(t,x)$.

We endow $\cS$ with the norm $\|F\|_{\bD^{1,2}}^{2}:=E|F|^2+E\|D F
\|_{\cH \cP}^{2}$. The operator $D$ can be extended to the space
$\bD^{1,2}$, the completion of $\cS$ with respect to $\|\cdot
\|_{\bD^{1,2}}$.

The following result is the analogue of Proposition 1.2.7 of \cite{nualart06}. %Its proof is in appendix.

\begin{proposition}
\label{Malliavin-rule-1} Let $F$ be a random variable given by
(\ref{Wiener-chaos-F}). If $F \in \bD^{1,2}$, then
$$D_{\bullet}F=\sum_{n \geq 1}n I_{n-1}(f_n(\cdot,\bullet)),$$
where $\bullet$ denotes the missing $(t,x)$ variable.
\end{proposition}

\noindent {\bf Proof:} It is enough to assume that $F=I_n(f_n)$, for
some symmetric elementary function $f_n$ of the form $f_n=\sum_{i_1,
\ldots,i_n=1}^{m}a_{i_1 \ldots i_n}1_{A_{i_1} \times \ldots \times
A_{i_n}}$, where $m \geq n$, $A_1, \ldots,A_m$ are pairwise-disjoint
bounded Borel sets in $\bR_{+} \times \bR^d$, and $a_{i_1 \ldots
i_n}=0$  if any two of the indices $i_1, \ldots,i_n$ are equal. Then
\begin{eqnarray}
\label{form-f-sym}f_n&=&\sum_{1 \leq i_1< \ldots<i_n \leq m}a_{i_1 \ldots i_n} \sum_{\rho \in S(\{i_1 \ldots i_n\})} 1_{A_{\rho(i_1)} \times \ldots \times A_{\rho(i_n)}} \\
\label{form-Inf}I_n(f_n)&=&n! \sum_{1 \leq i_1< \ldots <i_n \leq
m}a_{i_1 \ldots i_n}W(A_{i_1}) \ldots W(A_{i_n}),
\end{eqnarray}
where $S(\{i_1, \ldots,i_n\})$ denotes the set of all permutations
of $\{i_1, \ldots,i_n\}$ and $W(A)=W(1_{A})$. Using
(\ref{form-Inf}), the fact that $D(FG)=(DF)G+F(DG)$, and
$DW(\varphi)=\varphi$ for any $\varphi \in \cH \cP$, we infer that
$$D_{\bullet}F=n!  \sum_{1 \leq i_1< \ldots <i_n \leq m}a_{i_1 \ldots i_n} \sum_{j=1}^{n}1_{A_{i_j}}(\bullet)\left(\prod_{k \not = j}W(A_{i_k}) \right).$$

Using (\ref{form-f-sym}), we obtain:
\begin{eqnarray*}
I_{n-1}(f_n(\cdot, \bullet)) &=&
%\int_{(\bR_{+} \times \bR^d)^n}f_n(t_1,x_1, \ldots, t_{n-1},x_{n-1},\bullet)W(dt_1,dx_1) \ldots %W(dt_{n-1},dx_{n-1}) \\
%&=&
\sum_{1 \leq i_1< \ldots <i_n \leq m}a_{i_1 \ldots i_n} \sum_{\rho
\in S(\{i_1 \ldots i_n\})}
1_{A_{\rho(i_n)}}(\bullet) \prod_{j=1}^{n-1}W(A_{\rho(i_j)}) \\
&=& \sum_{1 \leq i_1< \ldots <i_n \leq m}a_{i_1 \ldots i_n}
\sum_{j=1}^{n} 1_{A_{i_j}}(\bullet) \sum_{\rho \in S(\{i_1 \ldots
i_n\} \verb2\2 \{i_j\})}
\prod_{k \not = j}W(A_{\rho(i_k)}) \\
&=& (n-1)!\sum_{1 \leq i_1< \ldots <i_n \leq m}a_{i_1 \ldots i_n}
\sum_{j=1}^{n} 1_{A_{i_j}}(\bullet) \prod_{k \not = j}W(A_{i_k}).
\end{eqnarray*}
$\Box$

The divergence operator $\delta$ is defined as the adjoint of the
operator $D$. The domain of $\delta$, denoted by $\mbox{Dom} \
\delta$, is the set of $u \in L^2(\Omega;\cH \cP)$ such that
$$|E \langle DF,u \rangle_{\cH \cP}| \leq c (E|F|^2)^{1/2}, \quad \forall F \in \bD^{1,2},$$
where $c$ is a constant depending on $u$. If $u \in {\rm Dom} \
\delta$, then $\delta(u)$ is the element of $L^2(\Omega)$
characterized by the following duality relation:
\begin{equation}
\label{duality} E(F \delta(u))=E\langle DF,u \rangle_{\cH \cP},
\quad \forall F \in \bD^{1,2}.
 \end{equation}

If $u \in \mbox{Dom} \ \delta$, we will use the notation
$$\delta(u)=\int_0^{\infty} \int_{\bR^d}u(t,x) W(\delta t, \delta x),$$
even if $u$ is not a function in $(t,x)$, and we say that
$\delta(u)$ is the Skorohod integral of $u$ with respect to $W$.

The next result %is the analogue of Proposition 1.3.7 of \cite{nualart06} and
gives an important calculus rule, which plays a crucial role in the
present article. This rule states, in particular, that the Skorohod
integral of a multiple Wiener integral of order $n$ coincides with a
multiple Wiener integral of order $n+1$, i.e.
$$\int_{\bR_{+} \times \bR^d} \left(\int_{(\bR_{+} \times \bR^d)^{n}} f_n(t_1,x_1, \ldots,t_n,x_n,t,x)W(dt_1,dx_1) \ldots W(dt_n, dx_n) \right)W(\delta t,\delta x)$$
\begin{equation}
\label{replace-d-by-delta}=\int_{(\bR_{+} \times \bR^d)^{n+1}}
f_n(t_1,x_1, \ldots,t_n,x_n,t,x)W(dt_1,dx_1) \ldots W(dt_n
dx_n)W(dt,dx).
\end{equation}

%The proof of Proposition \ref{Malliavin-rule-2} can be found in appendix.

\begin{proposition}
\label{Malliavin-rule-2} Assume that $u \in L^2(\Omega;\cH \cP)$ has
the Wiener chaos expansion:
\begin{equation}
\label{interpretation-rule2} u(\bullet)=\sum_{n \geq
0}I_n(f_n(\cdot,\bullet)),
\end{equation}
 where $\bullet$ denotes the missing $(t,x)$-variable, $\cdot$ denotes the missing $n$ variables $(t_1,x_1), \ldots, (t_n,x_n)$, and $f_n$ is symmetric and lies in  $\cH \cP^{\otimes n}$ (in the first $n$ variables).
 Then $u \in \mbox{Dom} \ \delta$ if and only if the series $\sum_{n \geq 0}I_{n+1}(\tilde f_n)$ converges in $L^2(\Omega)$, where $\tilde f_n$ is the symmetrization of $f_n$ in all $n+1$ variables. In this case,
 $$\delta(u)=\sum_{n \geq 0}I_{n+1}(\tilde f_n)=\sum_{n \geq 0}I_{n+1}(f_n).$$
 %In particular, %$I_n(f_n(\cdot,\bullet)) \in {\rm Dom} \ \delta$ and
% $\delta(I_n(f_n(\cdot,\bullet)))=I_{n+1}(f_n)$, i.e.
\end{proposition}

\begin{remark}
\label{interpretation-rule2-remark} {\rm (a) If $u(t,x)$ is a
function in $(t,x)$, relation (\ref{interpretation-rule2}) is
interpreted as follows: for any $(t,x) \in \bR_{+} \times \bR^d$,
\begin{equation}
\label{interpretation-rule2-function}u(t,x)=\sum_{n \geq 0}
I_n(f_{n}(\cdot,t,x)) \quad \mbox{in $L^2(\Omega)$}.
\end{equation}

\noindent (b) If $u(t,x)$ is a distribution in $(t,x)$, relation
(\ref{interpretation-rule2}) is interpreted as follows: for any
$\phi \in C_0^{\infty}(\bR_{+} \times \bR^d)$,
\begin{equation}
\label{interpretation-rule2-distribution}(u(\bullet),\phi)=\sum_{n
\geq 0} I_n((f_{n}(\cdot,\bullet), \phi) ) \quad \mbox{in
$L^2(\Omega)$}.
\end{equation}

\noindent (c) If $u(t,x)$ is a function in $t$ and a distribution in
$x$, relation (\ref{interpretation-rule2}) is interpreted as
follows: for any $t>0,\phi \in C_0^{\infty}(\bR^d)$,
$$(u(t,*),\phi)=\sum_{n \geq 0} I_n((f_{n}(\cdot,t,*), \phi) ) \quad
\mbox{in $L^2(\Omega)$},$$ where $*$ denotes the missing
$x$-variable. }
\end{remark}

\noindent {\bf Proof:} Using the same argument as in the proof of
Proposition 1.3.7 of \cite{nualart06}, it suffices to prove that for
any $G=I_n(g)$ with $g \in \cH \cP^{\otimes n}$ symmetric, we have:
\begin{equation}
\label{duality-G} E \langle DG, u \rangle_{\cH \cP}= E(I_{n}(\tilde
f_{n-1})G).
\end{equation}

Without loss of generality, we may assume that $g$ is a function in
all variables. By Proposition \ref{Malliavin-rule-1}, $DG$ is a
function given by
\begin{equation}
\label{Malliavin-derivative-G} D_{s,y}G=n I_{n-1} (g(\cdot,s,y))
\quad \forall s>0,y \in \bR^d.
\end{equation}

We consider separately the following three cases.

{\em Case 1.} $u(t,x)$ is a function in $(t,x)$. Using
(\ref{interpretation-rule2-function}),
(\ref{Malliavin-derivative-G}), the orthogonality of the Wiener
chaos spaces, and (\ref{isometry-In}), we obtain
\begin{eqnarray*}
E \langle DG, u \rangle_{\cH \cP} &=& \alpha_{H} E \int_{(\bR_{+} \times \bR^d)^2} u(t,x)(D_{s,y}G) |t-s|^{2H-2}f(x-y)dx dy dt ds \\
&=& n \alpha_{H}  \int_{(\bR_{+} \times \bR^d)^2} E(I_{n-1} (f_{n-1}(\cdot,t,x))I_{n-1} (g(\cdot,s,y))) \\ & & \quad \quad \quad |t-s|^{2H-2}f(x-y)dx dy dt ds \\
&=& n (n-1)! \ \alpha_{H}  \int_{(\bR_{+} \times \bR^d)^2} \langle f_{n-1}(\cdot,t,x),g(\cdot,s,y))) \rangle_{\cH \cP^{\otimes (n-1)}}\\
& & \quad \quad \quad |t-s|^{2H-2}f(x-y)dx dy dt ds \\
&=& n! \ \langle f_{n-1},g \rangle_{\cH \cP^{\otimes n}}=n! \
\langle \tilde f_{n-1},g \rangle_{\cH \cP^{\otimes n}} =E(I_n(\tilde
f_{n-1})I_n(g)),
%&=& I(I_n(\tilde f_{n-1})G).
\end{eqnarray*}
where for the second-last equality, we used the symmetry of $g$.
This proves (\ref{duality-G}).

{\em Case 2.} $u(t,x)$ is a distribution in $(t,x)$ (in
$\cS'(\bR^{d+1})$). In this case, we regularize $f_n$ as follows.
Let $\psi \in C_0^{\infty}(\bR^{d+1})$ be such that $\psi \geq 0$,
the support of $\psi$ is included in $(0,1) \times \{x \in \bR^d;|x|
\leq 1\}$ and $\int_{\bR^{d+1}}\psi(t,x)dt dx=1$. Let
$\psi_{\varepsilon}(t,x)=\varepsilon^{-d-1}
\psi(t/\varepsilon,x/\varepsilon)$ and
$f_{n,\varepsilon}(\cdot,\bullet):=
\psi_{\varepsilon}*f_n(\cdot,\bullet)$, where $*$ denotes the
convolution with respect to the missing $(t,x)$-variable, denoted by
$\bullet$. Note that $f_{n,\varepsilon}(\cdot,t,x)$ is a
%$C^{\infty}$
function in $(t,x)$ (see p. 245 of \cite{schwartz66}). Let
$$u_{\varepsilon}(t,x)= \sum_{n \geq
0}I_{n}(f_{n,\varepsilon}(\cdot,t,x)).$$

We claim that $u_{\varepsilon}=\psi_{\varepsilon}*u$. To see this,
note that for any $\phi \in C_0^{\infty}(\bR_{+} \times \bR^d)$,
\begin{eqnarray*}
((\psi_{\varepsilon}*u)(\bullet), \phi)&=& (u(\bullet),
\psi_{\varepsilon}* \tilde \phi)=\sum_{n \geq 0}
I_n((f_n(\cdot,\bullet),
\psi_{\varepsilon}*\tilde \phi)) \\
&= &\sum_{n \geq 0}I_n((f_{n,\varepsilon}(\cdot,\bullet), \phi))=
\sum_{n \geq 0}(I_n(f_{n,\varepsilon}(\cdot,\bullet)),
\phi)=(u_{\varepsilon}(\bullet),\phi),
% \quad \mbox{in $L^2(\Omega)$},
\end{eqnarray*}
where we used (\ref{interpretation-rule2-distribution}) for the
second equality, and the stochastic Fubini's theorem for the
second-last equality.

Applying the result of Case 1 to $u_{\varepsilon}$, we get:
\begin{equation}
\label{conv-DG-u} E \langle DG, u_{\varepsilon} \rangle_{\cH \cP}=
E(I_{n}(\tilde f_{n-1,\varepsilon})G).
\end{equation}
Relation (\ref{duality-G}) follows by letting $\varepsilon \to 0$.
On the left-hand side of (\ref{conv-DG-u}), we have:
\begin{eqnarray*}
E\|u_{\varepsilon}-u\|_{\cH \cP}^2 &=& a_{H,d} \int_{\bR}d\tau |\tau|^{1-2H} \int_{\bR^{d}}\mu(d\xi) |\cF u_{\varepsilon}(\tau,\xi)-\cF u(\tau,\xi)|^2 \\
&=& a_{H,d} \int_{\bR}d\tau |\tau|^{-(2H-1)} \int_{\bR^{d}}\mu(d\xi)
|\cF u(\tau,\xi)|^2 |\cF \psi_{\varepsilon}(\tau,\xi)-1|^2 \to 0,
\end{eqnarray*}
as $\varepsilon \to 0$, by the Dominated Convergence Theorem, where
$\cF$ denotes the Fourier transform in the $(t,x)$-variable and
$a_{H,d}$ is a constant depending on $H,d$ and $\mu$.

On the right-hand side of (\ref{conv-DG-u}), we have:
\begin{eqnarray*}\lefteqn{E|I_{n}(\tilde f_{n-1,\varepsilon}) -I_{n}(\tilde f_{n-1})|^2=n! \|\tilde f_{n-1,\varepsilon}-\tilde f_{n-1}\|_{\cH \cP^{\otimes n}}^2 =}\\
& & n! a_{H,d}^n \int_{\bR^n} d\tau_{1} \ldots d\tau_{n-1}d\tau
\prod_{i=1}^{n-1}|\tau_i|^{1-2H} |\tau|^{1-2H}
\int_{\bR^{nd}} \mu(d\xi_1) \ldots \mu(d\xi_{n-1})\mu(d\xi) \\
& & |\cF^{(n)} \tilde f_{n-1}(\tau_1,\xi_1, \ldots,\tau_{n-1},
\xi_{n-1},\tau,\xi)|^2 |\cF \psi_{\varepsilon}(\tau,\xi)-1|^2 \to 0,
\quad \mbox{as} \ \varepsilon \to 0,
\end{eqnarray*}
by the Dominated Convergence Theorem, where $\cF^{(n)}$ denotes the
Fourier transform in all $n$ variables $(t_1,x_1), \ldots,
(t_{n-1},x_{n-1}),(t,x)$.

{\em Case 3.} $u(t,x)$ is a function in $t$ and a distribution in
$x$. The argument is similar to Case 2, based on a regularization of
$u$ in space. We omit the details. $\Box$

We now return to our framework.

\begin{definition}
\label{definition-solution} We say that
 $u=\{u(t,x); t \geq 0, x \in \bR^d\}$ is a {\bf solution of (\ref{wave})} if $u(0,x)=1$ for all $x \in \bR^d$, and for any $t>0,x \in \bR^d$: \\
(i) $E|u(t,x)|^2<\infty$;\\
(ii) $u(t,x)$ is $\cF_t$-measurable, where $\cF_t=\sigma\{W_s(A); 0 \leq s \leq t,A \in \cB_{b}(\bR^d) \}$;\\
(iii) the process $v^{(t,x)}:=G(t-\cdot,x-*)u$ belongs to ${\rm Dom}
\ \delta$ and
$$u(t,x)=1+\delta(v^{(t,x)}).$$
Here, $\cdot$ denotes the missing $s$-variable, $*$ denotes the
missing $y$-variable, and $G(t-s,x-*)u(s, *)$ denotes the
multiplication of the distribution $G(t-s,x-*)$ with the function
$u(s,*)$, for any $s \in (0,t)$.
\end{definition}

The following result concludes our preliminary discussion.

\begin{theorem}
\label{prop-existence} Suppose that for any $t>0,x \in \bR^d, n \geq
1$, $f_{n}(\cdot,t,x) \in \cH \cP^{\otimes n}$, $f_{n}(\cdot,t,x)$
being the kernels introduced in Subsection \ref{def-kernels}. Then
equation (\ref{wave}) has a solution if and only if for any $t>0,x
\in \bR^d$,
\begin{equation}
\label{gen-nec-suf-cond} \mbox{the series $\sum_{n \geq
1}I_n(f_n(\cdot,t,x))$ converges in $L^2(\Omega)$}.
\end{equation} In this case,
the solution is given by: $u(0,x)=1$ for all $x \in \bR^d$, and
\begin{equation}
\label{def-sol} u(t,x)=1+\sum_{n \geq 1}I_n(f_n(\cdot,t,x)), \quad
\mbox{for all} \ t>0,x \in \bR^d.
\end{equation}
\end{theorem}

\noindent {\bf Proof:} %Suppose that (\ref{gen-nec-suf-cond}) holds,
%Let $u(t,x)$ be given by (\ref{def-sol}). It is clear that $u$
%satisfies properties (i)-(ii) in Definition
%\ref{definition-solution}. It remains to prove (iii).
Let $v^{(t,x)}$ be given by Definition \ref{definition-solution}. We
claim that $v^{(t,x)}$ has the Wiener chaos expansion:
\begin{equation}
\label{Wiener-exp-v} v^{(t,x)}(\bullet)=\sum_{n \geq
0}I_n(f_{n+1}(\cdot, \bullet,t,x)),
\end{equation} where $\bullet$
denotes the missing $(s,y)$-variable.

From (\ref{Wiener-exp-v}), by Proposition \ref{Malliavin-rule-2}, it
will follow that $v^{(t,x)} \in {\rm Dom} \ \delta$ if and only if
the series $\sum_{n \geq 0}I_{n+1}(f_{n+1}(\cdot,t,x))$ converges in
$L^2(\Omega)$, and in this case,
$$\delta(v^{(t,x)})=\sum_{n \geq 0}I_{n+1}(f_{n+1}(\cdot,t,x))=u(t,x)-1.$$

It remains to prove (\ref{Wiener-exp-v}). If $d \leq 2$, then
$G(t,x)$ is a function in $x$, and (\ref{Wiener-exp-v}) is clear,
since for any $s>0$ and $y \in \bR^d$,
\begin{eqnarray*}
v^{(t,x)}(s,y)&=&G(t-s,x-y)\sum_{n \geq
0}I_n(f_{n}(\cdot,s,y))=\sum_{n \geq
0}I_n(G(t-s,x-y)f_{n}(\cdot,s,y)) \\
&=& \sum_{n \geq 0}I_n(f_{n+1}(\cdot,s,y,t,x)).
\end{eqnarray*}

Suppose now that $d \geq 3$. Then $G(t,x)$ is a distribution in $x$.
Recalling the interpretation given to (\ref{Wiener-exp-v}) in Remark
\ref{interpretation-rule2-remark}.(c), we show that for any $s \in
(0,t)$ and $\phi \in C_0^{\infty}(\bR^d)$ fixed,
\begin{equation}
\label{Wiener-exp-v-phi} (v^{(t,x)}(s,*),\phi)=\sum_{n \geq 0}I_n(
(f_{n+1}(\cdot,s,*,t,x),\phi)) \quad \mbox{in $L^2(\Omega)$.}
\end{equation}

Since $v^{(t,x)}(s,*)$ is the product between the distribution
$G(t-s,x-*)$ and the function $u(s,*)$, the action of
$v^{(t,x)}(s,*)$ on $\phi$ is given by:
\begin{eqnarray*}
(v^{(t,x)}(s,*),\phi)&=&(G(t-s,x-*),\phi u(s,*))=[\phi
u(s,*)*G(t-s,*)](x) \\
&=& \sum_{n \geq 0}[\phi J_n(s,*)* G(t-s,*)](x),
\end{eqnarray*}
where $J_n(s,y)=I_n(f_n(\cdot,s,y))$, and we used (\ref{def-sol}).

To prove (\ref{Wiener-exp-v-phi}), it suffices to show that for any
$n \geq 0$,
\begin{equation}
\label{Wiener-exp-v-phi-step1} [\phi J_n(s,*)*
G(t-s,*)](x)=I_n((f_{n+1}(\cdot,s,*,t,x),\phi)).
\end{equation}

Let $G_{\varepsilon}$ be a regularization of $G$ in space, i.e.
$G_{\varepsilon}(s,*)=\psi_{\varepsilon}*G(s,*)$, where
$\psi_{\varepsilon}=\varepsilon^{-d} \psi(x/\varepsilon)$, $\psi \in
C_0^{\infty}(\bR^d)$, $\psi \geq 0$, ${\rm supp} \ \psi \subset \{x
\in \bR^d; |x| \leq 1\}$ and $\int_{\bR^d}\psi(x)dx=1$. Then $\cF
G_{\varepsilon}(s,*)(\xi)=\cF G(s,*)(\xi)\cF
\psi_{\varepsilon}(\xi)$.

Since $G_{\varepsilon}(s,y)$ is a function in $y$, for any $n \geq
0$,
\begin{equation}
\label{Wiener-exp-v-phi-step2} [\phi J_{n,\varepsilon}(s,*)*
G_{\varepsilon}(t-s,*)](x)=
I_n((f_{n+1,\varepsilon}(\cdot,s,*,t,x),\phi)),
\end{equation}
where $J_{n,\varepsilon}(s,y)=I_n(f_{n,\varepsilon}(\cdot,s,y))$,
$$f_{n,\varepsilon}(t_1,x_1,
\ldots,t_n,x_n,s,y)=G_{\varepsilon}(s-t_n,y-x_n) \ldots
G_{\varepsilon}(t_2-t_1,x_2-x_1)1_{\{t_1< \ldots<t_n<s\}},$$
$$f_{n+1,\varepsilon}(t_1,x_1, \ldots,t_n,x_n,s,y,t,x)=
G_{\varepsilon}(t-s,x-y)f_{n,\varepsilon}(t_1,x_1,
\ldots,t_n,x_n,s,y)1_{\{s<t\}}.$$

Relation (\ref{Wiener-exp-v-phi-step1}) follows by taking the limit
as $\varepsilon \to 0$ in (\ref{Wiener-exp-v-phi-step2}). This is
justified below.

On the left hand side of (\ref{Wiener-exp-v-phi-step2}), we use the
fact that for any $y \in \bR^d$,
\begin{equation}
\label{LHS-Wiener} E|J_{n,\varepsilon}(s,y)-J_n(s,y)|^2 \to 0, \quad
\mbox{as $\varepsilon \to 0$}.
\end{equation}
To see this, note that:
\begin{eqnarray*}
\lefteqn{E|J_{n,\varepsilon}(s,y)-J_n(s,y)|^2 =
E|I_n(f_{n,\varepsilon}(\cdot,s,y))-I_n(f_n(\cdot,s,y))|^2 = n!
\|(\tilde f_{n,\varepsilon}-\tilde f_n)(\cdot,s,y)\|_{\cH
\cP^{\otimes n}}^2 }
\\ & & = n! \alpha_H^n \int_{(0,s)^{2n}} \prod_{j=1}^{n}|t_j-s_j|^{2H-2}
\int_{\bR^{nd}} \cF (\tilde f_{n,\varepsilon}-\tilde
f_{n})(t_1,\cdot, \ldots, t_n,\cdot,s,y)(\bxi) \\
& & \quad \quad  \overline{\cF (\tilde f_{n,\varepsilon}-\tilde
f_{n})(s_1,\cdot, \ldots, s_n,\cdot,s,y)(\bxi)} \mu(d\xi_1) \ldots
\mu(d\xi_n) d{\bf t}d{\bf s},
\end{eqnarray*}
where $\bxi=(\xi_1,\ldots,\xi_n)$ and ${\bf t}=(t_1,\ldots,t_n)$.
Here, $\tilde f_n(\cdot,s,y)$ is the symmetrization of
$f_n(\cdot,s,y)$ in the first $n$ variables $(t_1,x_1), \ldots,
(t_n,x_n)$, and $\cF$ denotes the Fourier transform with respect to
the missing variables $x_1, \ldots,x_n$.

Using (\ref{Fourier-fn-formula}), one can prove that:
$$\cF (\tilde f_{n,\varepsilon}-\tilde
f_{n})(t_1,\cdot, \ldots, t_n,\cdot,s,y)(\bxi)=\cF \tilde
f_{n}(t_1,\cdot, \ldots, t_n,\cdot,s,y)(\bxi)[ \overline{\cF
\psi_{\varepsilon}(\xi_{\rho(1)})}
$$ $$\overline{\psi_{\varepsilon}(\xi_{\rho(1)}+\xi_{\rho(2)})}\ldots
\overline{\psi_{\varepsilon}(\xi_{\rho(1)}+\ldots+ \xi_{\rho(n)})}
-1],$$ where $\rho$ is a permutation such that $t_{\rho(1)}<
\ldots<t_{\rho(n)}$. Relation (\ref{LHS-Wiener}) follows by the
Dominated Convergence theorem, since $\cF \psi_{\varepsilon}(\xi)
\to 0$. The application of this theorem is justified since $|\cF
\psi_{\varepsilon}(\xi)| \leq 1$ and $\|\tilde f_n(\cdot,s,y)\|_{\cH
\cP^{\otimes n}}^2<\infty$.

On the right hand side of (\ref{Wiener-exp-v-phi-step2}), we use the
fact that:
\begin{equation}
\label{RHS-Wiener} E|I_n((f_{n+1,\varepsilon}(\cdot,s,*,t,x),\phi))-
I_n((f_{n+1}(\cdot,s,*,t,x),\phi))|^2 \to 0, \quad \mbox{as
$\varepsilon \to 0$}.
\end{equation}
To see this, note that
\begin{eqnarray}
\nonumber \lefteqn{E|I_n((f_{n+1,\varepsilon}(\cdot,s,*,t,x),\phi))-
I_n((f_{n+1}(\cdot,s,*,t,x),\phi))|^2 =n! \|
g_{\varepsilon}(\cdot,s,t,x)\|_{\cH \cP^{\otimes n}}^{2} } \\
\nonumber & & = n! \alpha_H^n \int_{(0,s)^{2n}}
\prod_{j=1}^{n}|t_j-s_j|^{2H-2} \int_{\bR^{nd}} \cF
g_{\varepsilon}(t_1,\cdot, \ldots,t_n,\cdot,s,t,x)(\bxi) \\
\label{RHS-Wiener-step1}& & \overline{\cF g_{\varepsilon}(s_1,\cdot,
\ldots,s_n,\cdot,s,t,x)(\bxi)} \mu(d\xi_1) \ldots \mu(d\xi_n)d{\bf
t} d{\bf s},
\end{eqnarray}
where $\bxi=(\xi_1, \ldots,\xi_{n})$, ${\bf t}=(t_1,\ldots,t_n)$ and
$$g_{\varepsilon}(\cdot,s,t,x):=((\tilde f_{n+1,\varepsilon}-\tilde
f_{n+1})(\cdot,s,*,t,x),\phi).$$ Here, the action of $\phi$ is on
the missing $y$-variable (denoted by $*$), $\tilde
f_{n+1}(\cdot,s,*,t,x)$ is the symmetrization of
$f_{n+1}(\cdot,s,*,t,x)$ in the first $n$ variables $(t_1,x_1),
\ldots, (t_n,x_n)$, and the Fourier transform is taken with respect
to the missing variables $x_1, \ldots,x_n$.

Note that
$$\cF g_{\varepsilon}(t_1,\cdot,
\ldots,t_n,\cdot,s,t,x)(\bxi)=\int_{\bR^d} \cF (\tilde
f_{n+1,\varepsilon}-\tilde f_{n+1})(t_1,\cdot,\ldots,
t_n,\cdot,s,*,t,x)(\bxi,\xi) \overline{\cF \phi(\xi)}d\xi,$$ where
the first Fourier transform under the integral is taken with respect
to the $n+1$ missing variables $x_1, \ldots,x_n,y$. Using
(\ref{Fourier-fn-formula}), one can prove that:
$$\cF (\tilde f_{n+1,\varepsilon}-\tilde f_{n+1})(t_1,\cdot,\ldots,
t_n,\cdot,s,*,t,x)(\bxi,\xi) = \cF \tilde f_n(t_1, \cdot, \ldots,
t_n,\cdot,s,x)(\bxi) k_{\varepsilon,\rho}(\bxi,\xi),$$ where $$
k_{\varepsilon,\rho}(\bxi,\xi)=e^{-i \xi \cdot x}\overline{\cF
G(t-s,\cdot)(\xi_1+ \ldots+\xi_{n}+\xi)} \
[\overline{\psi_{\varepsilon}(\xi_{\rho(1)})}
\ldots %\overline{\psi_{\varepsilon}(\xi_{\rho(1)}+ \ldots +\xi_{\rho(n)})}\
\overline{\psi_{\varepsilon}(\xi_{\rho(1)}+ \ldots
+\xi_{\rho(n)}+\xi)}- 1],$$ and $\rho$ is the permutation for which
$t_{\rho(1)}< \ldots <t_{\rho(n)}$. Hence
\begin{eqnarray}
\nonumber \cF g_{\varepsilon}(t_1,\cdot,
\ldots,t_n,\cdot,s,t,x)(\bxi)&=&\cF \tilde f_n(t_1, \cdot, \ldots,
t_n,\cdot,s,x)(\bxi)\int_{\bR^d}k_{\varepsilon,\rho}(\bxi,\xi)
\overline{\cF \phi (\xi)}d\xi \\
\label{RHS-Wiener-step2} &=:& \cF \tilde f_n(t_1, \cdot, \ldots,
t_n,\cdot,s,x)(\bxi)K_{\varepsilon,\rho}(\bxi).
\end{eqnarray}

Since $\cF \psi_{\varepsilon}(\xi) \to 0$, by the Dominated
Convergence Theorem, it follows that $K_{\varepsilon,\rho}(\bxi) \to
0$ as $\varepsilon \to 0$. (To justify this, we use
(\ref{FG-bounded}) and $|\cF \psi_{\varepsilon}(\xi)| \leq 1$.)

Relation (\ref{RHS-Wiener}) follows from (\ref{RHS-Wiener-step1})
and (\ref{RHS-Wiener-step2}), again by the Dominated Convergence
Theorem, whose application is justified by the fact that
$|K_{\varepsilon,\rho}(\bxi)| \leq 2 C_{t} \int_{\bR^d}|\cF
\phi(\xi)|d\xi=:C_{t,\phi}$, and $\|\tilde f_{n}(\cdot,s,x)\|_{\cH
\cP^{\otimes n}}^2<\infty$. $\Box$

\begin{remark}
{\rm Let $u_0(t,x)=1$ and
$u_{n}(t,x)=1+\sum_{k=1}^{n}I_k(f_k(\cdot,t,x))$ for $n \geq 1$.
Let $v_n^{(t,x)}=G(t-\cdot,x-*)u_n$ for any $n \geq 0$. Using the same
argument as above, %(based on Proposition \ref{Malliavin-rule-2}),
one can show that $\delta(v_n^{(t,x)})=u_{n+1}(t,x)-1$, i.e
$$u_{n+1}(t,x)=1+\int_0^t \int_{\bR^d}G(t-s,x-y)u_{n}(s,y)W(\delta s, \delta y), \quad \forall n \geq 0.$$
In other words, $\{u_n\}_{n \geq 0}$ plays the role of the Picard's
iteration sequence used in the case $H=1/2$.
 }
\end{remark}

\section{Existence of the Solution}
\label{existence}

In this section, we examine condition (\ref{gen-nec-suf-cond}), in
the particular case when $f$ is the Riesz kernel of order $\alpha
\in (0,d)$, i.e. $f(x)=c_{\alpha,d}|x|^{-(d-\alpha)}$,
$\mu(d\xi)=|\xi|^{-\alpha}d\xi$.

Our main result shows that $\alpha>d-2$ is a sufficient for
(\ref{gen-nec-suf-cond}), and hence a sufficient condition for the
existence of a solution to (\ref{wave}) (by Theorem
\ref{prop-existence}).

Note that $\alpha>d-2$ is also the sufficient condition for the
existence of a solution to (\ref{wave}), in the case when $H=1/2$
and the solution is interpreted using a martingale measure
stochastic integral (see Theorem 5.1 of \cite{conus-dalang09}).
%For a general statement of this result,

Due to the orthogonality of the Wiener chaos spaces and
(\ref{isometry-In}), condition (\ref{gen-nec-suf-cond}) is
equivalent to:
\begin{equation}
\label{S-finite} S(t):=1+\sum_{n \geq 1}n! \|\tilde
f_n(\cdot,t,x)\|_{\cH \cP^{\otimes n}}^2<\infty,
\end{equation}
where $\tilde f_n(\cdot,t,x)$ is the symmetrization of
$f_n(\cdot,t,x)$ in the first $n$ variables $(t_1,x_1), \linebreak
\ldots, (t_n,x_n)$. In this case, $E|u(t,x)|^2=S(t)$.

We begin with the calculation of $\|\tilde f_n(\cdot,t,x)\|_{\cH
\cP^{\otimes n}}^2$.
 At the same time, this calculation will show that $f_n(\cdot,t,x) \in \cH \cP^{\otimes n}$ (see Remark \ref{rem-varphi-in-HP}).
 (By abuse of notation, we use $\|\cdot \|_{\cH \cP^{\otimes n}}$, even if we do not know yet that
$\tilde f_n(\cdot,t,x) \in \cH \cP^{\otimes n}$.)

Note that
$$\|\tilde f_n(\cdot,t,x)\|_{\cH
\cP^{\otimes n}}^2=\alpha_H^n \int_{(0,t)^{2n}}
\prod_{j=1}^{n}|t_j-s_j|^{2H-2} \int_{\bR^d} \cF \tilde
f_n(t_1,\cdot, \ldots,t_n, \cdot,t,x)(\bxi) $$ $$\overline{\cF
\tilde f_n(s_1,\cdot, \ldots,s_n, \cdot,t,x)(\bxi)} \mu(d\xi_1)
\ldots \mu(d\xi_n) d{\bf t} d{\bf s},$$ where
$\bxi=(\xi_1,\ldots,\xi_n)$, ${\bf t}=(t_1, \ldots, t_n)$ and ${\bf
s}=(s_1,\ldots,s_n)$. Let
$$g_{\bf t}^{(n)}(\cdot,t,x):=n! \tilde f_n(t_1, \cdot, \ldots,t_n,\cdot,t,x).$$
Hence,
$$\cF \tilde f_n(t_1, \cdot, \ldots,t_n,\cdot,t,x)(\bxi)=\frac{1}{n!}\cF g_{\bf
t}^{(n)}(\cdot,t,x)(\bxi),$$ and
$$ \|\tilde f_n(\cdot \ , t,x)\|_{\cH \cP^{\otimes
n}}^{2}=\frac{1}{(n!)^2}\alpha_H^n \int_{[0,t]^{2n}}
\prod_{j=1}^{n}|s_j-t_j|^{2H-2} \tilde \psi^{(n)}({\bf t}, {\bf s})
d{\bf t} d {\bf s},$$ where
\begin{equation}
\label{def-psi*-n} \tilde \psi^{(n)}({\bf t}, {\bf s}) :=
%\langle g_{\bf s}^{(n)} , g_{\bf t}^{(n)} \rangle_{\cP(\bR^d)^{\otimes n}}
\int_{\bR^{nd}} \cF g_{\bf t}^{(n)}(\cdot,t,x) (\bxi) \overline{\cF
g_{\bf s}^{(n)}(\cdot,t,x)(\bxi)} \mu(d\xi_1) \ldots \mu(d\xi_n).
\end{equation}

We let
$$\tilde \alpha_n(t):= (n!)^2 \|\tilde f_n(\cdot \ , t,x)\|_{\cH \cP^{\otimes
n}}^{2}.$$ With this notation, relation (\ref{S-finite}) becomes:
\begin{equation}
\label{S-finite-2} S(t)=\sum_{n \geq 0}\frac{1}{n!}\tilde
\alpha_n(t)<\infty.
\end{equation}

We proceed to the evaluation of $\tilde \alpha_n(t)$, which relies
on the evaluation of $\tilde \psi^{(n)}({\bf t},{\bf s})$. Using
relation (\ref{Fourier-fn-formula}), one can prove that:
 $$\cF g_{\bf t}^{(n)}(\cdot,t,x)(\bxi)=e^{-i
(\xi_1+\ldots +\xi_n)\cdot x} \sum_{\rho \in S_n} \overline{\cF
G(u_1,\cdot \ )(\xi_{\rho(1)})} \ \overline{\cF G(u_2,\cdot \
)(\xi_{\rho(1)}+\xi_{\rho(2)})} $$
\begin{equation}
\label{Fourier-transform-gt} \ldots \overline{\cF G(u_n,\cdot \
)(\xi_{\rho(1)}+\ldots+\xi_{\rho(n)})} 1_{\{t_{\rho(1)}<
\ldots<t_{\rho(n)} \}},
\end{equation}
where $S_n$ is the set of all permutations of $\{1, \ldots,n\}$ and
$u_j=t_{\rho(j+1)}-t_{\rho(j)}$ for $1 \leq j \leq n$, with
$t_{\rho(n+1)}=t$.

In the argument below, since there is no risk of confusion, we omit
writing the variable $(t,x)$ of $g_{\bf t}^{(n)}(\cdot,t,x)$.

It is known that, for any $t>0$, $G(t,\cdot)$ is a distribution with
rapid decrease in $\cS'(\bR^d)$, whose Fourier transform is given
by: (see e.g. \cite{taylor96})
$$\cF G(t,\cdot)(\xi)=\frac{\sin (t|\xi|)}{|\xi|}, \quad \forall \xi \in \bR^d.$$

The following central result was announced in the introduction, and
will allow us to estimate $\tilde \psi^{(n)}({\bf t}, {\bf s})$.

\begin{lemma}
\label{hu-lemma}Assume that $d-2<\alpha<d$. Then,
$$I:=\int_{\bR^d}\frac{\sin^2(t|\xi|)}{|\xi|^2}|\xi-\eta|^{-\alpha}d\xi\leq
C_{\alpha,d}t^{\alpha-d+2}, \quad \mbox{for any} \ t>0,\eta \in
\bR^d.$$
\end{lemma}

\noindent {\bf Proof:} Using the change of variable $\xi'=t\xi$, we
obtain,
$$I=t^{\alpha-d+2} \int_{\bR^d}
\frac{\sin^2(|\xi'|)}{|\xi'|^2}|\xi'-t\eta|^{-\alpha}d\xi'.$$

We claim that:
$$I(a):= \int_{\bR^d}
\frac{\sin^2(|\xi|)}{|\xi|^2}|a-\xi|^{-\alpha}d\xi \leq
C_{\alpha,d},\quad \forall a \in \bR^d.$$ To see this, we change the
variable $a-\xi$ into $\xi$, and we write
$$I(a)=\int_{|\xi| \leq 1}
\frac{\sin^2(|\xi-a|)}{|\xi-a|^2}|\xi|^{-\alpha}d\xi +\int_{|\xi|>1}
\frac{\sin^2(|\xi-a|)}{|\xi-a|^2}|\xi|^{-\alpha}d\xi=:I_1(a)+I_2(a).$$

For $I_1(a)$, we use the fact that $|\frac{\sin x}{x}| \leq 1$ for
any $x>0$. Hence
$$I_1(a) \leq  \int_{|\xi| \leq 1}|\xi|^{-\alpha}d\xi=c_d \int_0^1
\lambda^{-\alpha+d-1}d\lambda=c_d \frac{1}{d-\alpha}.$$

For $I_2(a)$, we use the fact that
$$\frac{\sin^2 (t|\xi|)}{|\xi|^2} \leq 2(t^2+1)\frac{1}{1+|\xi|^2}, \quad \forall t>0, \forall \xi \in \bR^d.$$
(see p. 81 of \cite{sanzsole05}). In our case, $t=1$. Hence
$$I_2(a) \leq 4 \int_{|\xi|>1}
\frac{1}{1+|\xi-a|^2}|\xi|^{-\alpha}d\xi\leq 4 \sup_{a \in \bR^d}
\int_{\bR^d}\frac{1}{1+|\xi-a|^2}|\xi|^{-\alpha}d\xi.$$

Finally, we observe that $\alpha>d-2$ is equivalent to
$\int_{\bR^d}\frac{1}{1+|\xi|^2}|\xi|^{-\alpha}d\xi<\infty$, which
in turn is equivalent to $$\sup_{a \in \bR^d}
\int_{\bR^d}\frac{1}{1+|\xi-a|^2}|\xi|^{-\alpha}d\xi<\infty$$ (see
(5.5) of \cite{conus-dalang09}). $\Box$

\vspace{3mm}

Based on the previous lemma, we estimate $\tilde \psi^{(n)}({\bf
t},{\bf s})$.

%The next result is the analogue of Lemma 3.4 of \cite{balan-tudor09}
%in the case of the wave equation.

\begin{lemma}
\label{bound-psi-n} If $f$ is the Riesz kernel of order
$\alpha>d-2$, then for any ${\bf t}, {\bf s} \in [0,t]^n$,
$$\tilde \psi^{(n)}({\bf t}, {\bf s}) \leq C_{\alpha,d}^{n}[\beta({\bf
t}) \beta({\bf s})]^{(\alpha-d+2)/2},$$ where $\beta({\bf
t})=\prod_{j=1}^{n}(t_{\rho(j+1)}-t_{\rho(j)})$, $\beta({\bf
s})=\prod_{j=1}^{n}(s_{\sigma(j+1)}-s_{\sigma(j)})$, and the
permutations $\rho$ and $\sigma$ of $\{1,\ldots,n\}$ are chosen such
that
\begin{equation}
\label{permut-rho-sigma} t_{\rho(1)}<t_{\rho(2)}< \ldots <
t_{\rho(n)} \quad \mbox{and} \quad s_{\sigma(1)}<s_{\sigma(2)}<
\ldots < s_{\sigma(n)},
\end{equation} with
$t_{\rho(n+1)}=s_{\sigma(n+1)}=t$.
\end{lemma}

\noindent {\bf Proof:} By the Cauchy-Schwartz inequality,
%\begin{equation}
%\label{CS-psi}
$$\tilde \psi^{(n)}({\bf t}, {\bf s}) \leq \tilde \psi^{(n)}({\bf
t}, {\bf t})^{1/2}\tilde \psi^{(n)}({\bf s}, {\bf s})^{1/2}.$$
%\end{equation}

Let $u_j=t_{\rho(j+1)}-t_{\rho(j)}$ for $j=1, \ldots, n$. Using
(\ref{def-psi*-n}) and (\ref{Fourier-transform-gt}), we obtain:
\begin{eqnarray*}
\tilde \psi^{(n)}({\bf t}, {\bf t}) &=& \int_{\bR^{nd}} |\cF
g_{\bf t}^{(n)}(\bxi)|^2 \mu(d\xi_1) \ldots \mu(d\xi_n) \\
&=& \int_{\bR^{nd}}|\cF G(u_1,\cdot \ )(\xi_{\rho(1)})|^2 \ldots
|\cF G(u_n,\cdot \ )(\xi_{\rho(1)}+\ldots+\xi_{\rho(n)})|^2
\mu(d\xi_1) \ldots \mu(d\xi_n) \\
&=& \int_{\bR^{nd}}|\cF G(u_1,\cdot \ )(\xi_{1}')|^2 \ldots |\cF
G(u_n,\cdot \ )(\xi_{1}'+\ldots+\xi_{n}')|^2 \mu(d\xi_1') \ldots
\mu(d\xi_n') \\
&=& \int_{\bR^{nd}} \frac{\sin^2(u_1|\xi_1'|)}{|\xi_1'|^2}\cdot
\frac{\sin^2(u_2|\xi_1'+\xi_2'|)}{|\xi_1'+\xi_2'|^2} \ldots
\frac{\sin^2(u_n|\xi_1'+\ldots
+\xi_n'|)}{|\xi_1'+\ldots+\xi_n'|^2}|\xi_1'|^{-\alpha}\ldots|\xi_n'|^{-\alpha}
d\xi_1'\ldots d\xi_n',
\end{eqnarray*}
where we used the change of variable $\xi_j'=\xi_{\rho(j)},
j=1,\ldots,n$.

We now use the change of variable
$$\eta_j=\xi_1'+\ldots+\xi_j', \quad j=1,\ldots,n.$$
The inverse transformation is: $\xi_1'=\eta_1,
\xi_j'=\eta_j-\eta_{j-1},j=2,\ldots,n$. We get
\begin{eqnarray*}
\psi^{*(n)}({\bf t}, {\bf t}) &=& \int_{\bR^{d}}d\eta_1
\frac{\sin^2(u_1|\eta_1|)}{|\eta_1|^2}|\eta_1|^{-\alpha}\int_{\bR^d}d\eta_2
\frac{\sin^2(u_2|\eta_2|)}{|\eta_2|^2}|\eta_2-\eta_1|^{-\alpha}
\ldots \\
& & \int_{\bR^d}d\eta_n\frac{\sin^2(u_n|\eta_n|)}{|\eta_n|^2}
|\eta_n-\eta_{n-1}|^{-\alpha}.
\end{eqnarray*}

Using Lemma \ref{hu-lemma} iteratively, we obtain $\tilde
\psi^{(n)}({\bf t}, {\bf t}) \leq C_{\alpha,d}^n (u_1 \ldots
u_n)^{\alpha-d+2}$. The result follows. $\Box$ %by (\ref{CS-psi}).

%The next result is the analogue of Proposition 3.6 in the case of
%the wave equation.

\begin{proposition}
\label{bound-alpha-n}  If $f$ is the Riesz kernel of order
$\alpha>d-2$, then for any $t>0$ and $n \geq 1$,
\begin{equation}
\label{bound-alpha-n-wave} \tilde \alpha_n(t) \leq C(t)^n
\frac{1}{(n!)^{\alpha-d+2}},
\end{equation}
 where
$C(t)=C_{\alpha,d,H}t^{2H+\alpha-d+2}$.
\end{proposition}

\noindent {\bf Proof:} Let $h=(\alpha-d+2)/(2H)$. As in the proof of
Proposition 3.6 of \cite{balan-tudor09}, using Lemma
\ref{bound-psi-n} and inequality (16) of \cite{balan-tudor09}, we
obtain that:
\begin{eqnarray*}
\tilde \alpha_n(t) & \leq & C_{\alpha,d,H}^n
(n!)^{2H}\left(\int_{0<s_1<s_2< \ldots <s_n<t}
[(t-s_n) \ldots (s_2-s_1)]^h d{\bf s}\right)^{2H} \\
& \leq & C_{\alpha,d,H}^n t^{n(1+h)2H}
\left(\frac{n!}{\Gamma(n(1+h)+1)} \right)^{2H} \\
& \leq & C_{\alpha,d,H}^n t^{n(1+h)2H} \frac{1}{(n!)^{\alpha-d+2}}.
\end{eqnarray*}

For the last inequality above, we used the fact that for $a>0$,
$\Gamma(an+1)=C_n (n!)^{a}$, where $C_n$ is a constant such that
$\lambda^{-n} \leq C_n \leq \lambda^n$ for some $\lambda>1$.
% (see p. 236-237 of \cite{hu01}).
$\Box$

The existence of the solution is immediate.

\begin{proposition}
If $f$ is the Riesz kernel of order $\alpha>d-2$, then condition
(\ref{gen-nec-suf-cond}) holds, and consequently, equation
(\ref{wave}) has a solution.
\end{proposition}

\noindent {\bf Proof:} As we mentioned earlier, condition
(\ref{gen-nec-suf-cond}) is equivalent to (\ref{S-finite-2}), which
in turn is satisfied, since by Proposition \ref{bound-alpha-n},
$$S(t)=\sum_{n \geq 0}\frac{1}{n!}\tilde \alpha_n(t) \leq \sum_{n \geq
0}\frac{C(t)^n}{(n!)^{\alpha-d+3}}<\infty.$$ The second statement
follows by Theorem \ref{prop-existence}. $\Box$

\vspace{3mm}

Since $C(t)$ is an increasing function in $t$, the previous argument
shows that:
\begin{equation}
\label{S-finite-2-strong} S(t)=\sum_{n \geq 0}\frac{1}{n!}\tilde
\alpha_n(t) \leq C_T<\infty, \quad \forall t \in [0,T],
\end{equation}
for any $T>0$, i.e. $\sup_{(t,x) \in [0,T] \times
\bR^d}E|u(t,x)|^2<\infty$ for all $T>0$.

\begin{remark}
\label{remark-heat} {\rm In the case of the heat equation, it was
shown in \cite{balan-tudor09} that, if $f$ is the Riesz kernel of
order $\alpha>d-4H$, then for any $t>0$ and $n\geq 1$,
%\begin{equation}
%\label{bound-alpha-n-heat}
$$\tilde \alpha_n(t) \leq C(t)^n
\frac{1}{(n!)^{-(d-\alpha)/2}},$$
%\end{equation}
where $C(t)=C_{\alpha,d,H} t^{2H-(d-\alpha)/2}$. In this case,
(\ref{S-finite-2}) holds % $\sum_{n \geq 0} \frac{1}{n!}\alpha_n(t)<\infty$
if $\alpha>d-2$.}
\end{remark}

\begin{remark}
\label{remark-alpha-non-sym} {\rm Using the same method as above,
one can prove that $\sum_{n \geq 0}\frac{1}{n!}\alpha_n(t) \leq C_T<
\infty$ for all $t \in [0,T]$, where
$\alpha_n(t)=(n!)^2\|f_n(\cdot,t,x)\|_{\cH \cP^{\otimes n}}^{2}$. }
\end{remark}

\section{Moments of the Solution}
\label{moments}

In this section, we show that the solution is
$L^2(\Omega)$-continuous and has uniformly bounded moments of order
$p \geq 1$. With the obvious modifications, the results presented in
this section remain valid for the heat equation (see Remark
\ref{remark-heat}).

Let $u(t,x)=\sum_{n \geq 0}J_n(t,x),$ where $J_n(t,x)$ is the
projection of $u(t,x)$ on the Wiener chaos $\cH \cP_{n}$. By the
orthogonality of the $J_n(t,x)$'s, we have:
\begin{equation}
\label{moment2-u}
 E|u(t,x)|^2=\sum_{n \geq 0}E|J_n(t,x)|^2.
\end{equation}
Note that
\begin{equation}
 \label{moment2-Jn} E|J_n(t,x)|^2= E|I_n(\tilde f_n(\cdot \ ,
t,x))|^2=%\|f_n(\cdot \ , t,x)\|_{\cH\cP^{\otimes n}}^{2}=
 \frac{1}{n!}\tilde \alpha_n(t).
 \end{equation}

It is known that, %if $W=\{W(h);h \in \cH\}$ is an isonormal Gaussian process, then
for any $1<p<q<\infty$, the norms $\|\cdot \|_{p}$ and $\|\cdot
\|_{q}$ are equivalent on any Wiener chaos $\cH \cP_n$, where
$\|\cdot\|_p$ denotes the norm in $L^p(\Omega)$. This is a
consequence of the hypercontractivity property of the
Ornstein-Uhlenbeck semigroup $(T_t)_{t \geq 0}$, defined by: $$T_t
F=\sum_{n \geq 0}e^{-nt}J_n F, \quad F \in L^2(\Omega),$$ where we
denote by $J_n F$ the projection of $F$ on the $n$-th Wiener chaos
$\cH \cP_n$. The property says that for any $p>1$ and $t>0$,
$$\|T_t F\|_{q(t)} \leq \|F \|_{p},$$
where $q(t)=e^{2t}(p-1)+1$ (see Theorem 1.4.1 of \cite{nualart06}).
Hence, for any $1<p<q<\infty$ and for any $F \in \cH \cP_n$,
$$e^{-nt}\|F\|_{q}=\|T_t F \|_{q} \leq \|F\|_p,$$ where $t>0$ is
chosen such that $q=e^{2t}(p-1)+1$. In particular, for any $p>2$ and
for any $F \in \cH \cP_n$,
\begin{equation}
\label{norms-p-q-equiv}\|F\|_{p} \leq  e^{nt}\|F\|_2=(p-1)^{n/2}
\|F\|_{2},
\end{equation}
where $t>0$ is chosen such that $p=e^{2t}+1$.

Applying these results in our case, we obtain the following result.

\begin{theorem}
Let $f$ be the Riesz kernel of order $\alpha>d-2$ and $u$
be the solution of (\ref{wave}). % the wave (or heat)  equation.
Then $u$ is $L^2(\Omega)$-continuous, and for any $p \geq 1$, $T>0$
\begin{equation}
\label{Lp-norm-bded} \sup_{t \leq T} \sup_{x \in \bR^d}E|u(t,x|^p
<\infty.
\end{equation}
\end{theorem}

\noindent {\bf Proof:} We apply (\ref{norms-p-q-equiv}) for
$F=J_n(t,x) \in
\cH \cP_n$. Using (\ref{moment2-Jn}), we obtain that: %for any $t>0,x \in \bR^d$,
$$\|J_n(t,x)\|_p \leq (p-1)^{n/2}\|J_n(t,x)\|_2=(p-1)^{n/2}\left(\frac{1}{n!}
\tilde \alpha_n(t)\right)^{1/2}.$$

%In the case of the wave equation,
Using (\ref{bound-alpha-n-wave}), we obtain that:
$$\sum_{n \geq 0}\|J_n(t,x)\|_{p} \leq \sum_{n \geq 0} (p-1)^{n/2}
\left\{ C(t)^n \frac{1}{(n!)^{\alpha-d+3}}\right\}^{1/2}<\infty.$$
%whereas in the case of the heat equation, using
%(\ref{bound-alpha-n-heat}), we obtain that:
%$$\sum_{n \geq 0}\|J_n(t,x)\|_{p} \leq \sum_{n \geq 0} (p-1)^{n/2}
%\left\{ C(t)^n
%\frac{1}{(n!)^{1-(d-\alpha)/2}}\right\}^{1/2}<\infty.$$

%In both cases,
Since $\alpha(t)$ does not depend on $x$ and $C(t)$ is an increasing
function of $t$, we have: for any $T>0$,
\begin{equation}
\label{sup-Jn-bded} \sum_{n \geq 0} \sup_{t \leq T} \sup_{x \in
\bR^d} \|J_n(t,x)\|_{p}\leq C_{T,p}<\infty.
\end{equation}

From here, we conclude that for any $(t,x) \in \bR_{+} \times
\bR^d$, the sequence $\{u_n(t,x)=\sum_{k=0}^{n}J_k(t,x), n \geq 0\}$
is Cauchy in $L^p(\Omega)$, since
$$\|u_n(t,x)-u_m(t,x)\|_p \leq \sum_{k=m+1}^{n}\|J_k(t,x)\|_{p}
\to 0, \quad \mbox{as} \ n,m \to \infty,n>m.$$ Therefore, there
exists a random variable $v(t,x) \in L^p(\Omega)$ such that
$u_n(t,x) \to v(t,x)$ in $L^p(\Omega)$. But $u_{n}(t,x) \to u(t,x)$
in $L^2(\Omega)$, and hence $u(t,x)=v(t,x)$ a.s. Using
(\ref{sup-Jn-bded}), one can show that
\begin{equation}
\label{unif-cont} u_n(t,x) \to u(t,x) \ \mbox{in $L^p(\Omega)$,
uniformly in $(t,x) \in [0,T] \times \bR^d$}
\end{equation}
%since
%$$\|u_n(t,x)-u(t,x)\|_p =\|\sum_{k \geq n+1} J_k(t,x)\|_p \leq
%\sum_{k \geq n+1}\|J_k(t,x)\|_p<\varepsilon.$$
and $\|u_n(t,x)\|_p \leq \sum_{k=0}^{n}\|J_k(t,x)\|_p \leq C_{T,p}$
for all $(t,x) \in [0,T] \times \bR^d, n \geq 0$. Taking $n \to
\infty$, we obtain (\ref{Lp-norm-bded}).

By Lemma \ref{Jn-cont} below, $J_n$ is $L^2(\Omega)$-continuous.
Hence $u_n$ is $L^2(\Omega)$-continuous. Due to (\ref{unif-cont}),
it follows that $u$ is $L^2(\Omega)$-continuous. $\Box$

\begin{lemma}
\label{Jn-cont}a) For any $n \geq 1$ and $t>0$,
$$E|J_n(t+h,x)-J_n(t,x)|^2 \to 0 \ \mbox{as $h \to 0$, uniformly in $x
\in \bR^d$}.$$

b) For any $n \geq 1, t>0,x \in \bR^d$,
$$E|J_n(t,x+z)-J_n(t,x)|^2 \to 0 \ \mbox{as $|z| \to 0,z \in \bR^d$}.$$
\end{lemma}

\noindent {\bf Proof:} %We treat only the wave equation, the heat equation being similar.
%The techniques of the proof are similar to the white noise case.
a) Suppose that $h \in [0,1]$. (The case $h<0$ is similar.) Then,
\begin{eqnarray*}
E|J_n(t+h,x)-J_n(t,x)|^2 &=& E|I_n(\tilde f_n(\cdot  , t+h,x)-\tilde
f_n(\cdot
,t,x))|^2 \\
&=& n! \ \|\tilde f_n(\cdot ,t+h,x)-\tilde f_n(\cdot ,t,x)\|_{\cH
\cP^{\otimes n}}^{2} \\
& \leq & \frac{2}{n!}(E_1(t,h)+E_2(t,h)),
\end{eqnarray*}
where
\begin{eqnarray}
\label{def-E1} E_1(t,h)&:=& (n!)^2\|\tilde f_n(\cdot ,
t+h,x)1_{[0,t]^n}-\tilde f_n(\cdot
,t,x)\|_{\cH \cP^{\otimes n}}^{2} \\
\label{def-E2} E_2(t,h) &:=& (n!)^2\|\tilde f_n(\cdot
,t+h,x)1_{[0,t+h]^n \verb2\2 [0,t]^n}\|_{\cH \cP^{\otimes n}}^{2}.
\end{eqnarray}

We treat $E_1(t,h)$ first. Note that
\begin{equation}
\label{formula-E1} E_1(t,h) =\alpha_H^n \int_{[0,t]^{2n}}
\prod_{j=1}^{n}|t_j-s_j|^{2H-2} \psi_{h}^{(n)}({\bf t}, {\bf
s})d{\bf t} d{\bf s},
\end{equation} where
$$\psi_{h}^{(n)}({\bf t},{\bf s}) =\int_{\bR^{nd}}
\cF(g_{\bf t}^{(n)}(\cdot,t,x+h)-g_{\bf
t}^{(n)}(\cdot,t,x))(\bxi)\overline{\cF(g_{\bf
s}^{(n)}(\cdot,t,x+h)-g_{\bf s}^{(n)}(\cdot,t,x))(\bxi)}\mu(d\xi_1)
\ldots \mu(d\xi_n)$$
%\langle g_{\bf t}^{(n)}(\cdot \ ; t+h,x)-g_{\bf t}^{(n)}(\cdot \ ; t,x) , \ g_{\bf s}^{(n)}(\cdot \ ;
%t+h,x)-g_{\bf s}^{(n)}(\cdot \ ; t,x)\rangle_{\cP(\bR^d)^{\otimesn}}.$$

\noindent By the Cauchy-Schwartz inequality, $\psi_{h}^{(n)}({\bf
t},{\bf s}) \leq \psi_{h}^{(n)}({\bf t},{\bf t})^{1/2} \cdot
\psi_{h}^{(n)}({\bf s},{\bf s})^{1/2}$.

To evaluate $\psi_{h}^{(n)}({\bf t},{\bf t})$, we use
(\ref{Fourier-transform-gt}), denoting
$u_j=t_{\rho(j+1)}-t_{\rho(j)}$, when $0<t_{\rho(1)}< \ldots
<t_{\rho(n)}<t_{\rho(n+1)}=t$:
\begin{eqnarray*}
\psi_{h}^{(n)}({\bf t},{\bf t})
%\|g_{\bf t}^{(n)}(\cdot \ ; t+h,x)-g_{\bf t}^{(n)}(\cdot \ ; t,x)\|_{\cP(\bR^d)^{\otimes n}}^{2} \\
&=& \int_{\bR^{nd}} |\cF (g_{\bf t}^{(n)}(\cdot \ , t+h,x) - g_{\bf
t}^{(n)}(\cdot \ , t,x))(\bxi)|^2 \mu(d\xi_1) \ldots
\mu(d\xi_n) \\
&=& \int_{\bR^{nd}} |\cF G(u_1, \cdot \ )(\xi_{\rho(1)})|^2 \ldots
|\cF G(u_{n-1}, \cdot \ )(\xi_{\rho(1)}+ \ldots +\xi_{\rho(n-1)})|^2
\\
& & |\cF [G(u_n+h, \cdot \ )-G(u_n , \cdot \ )](\xi_{\rho(1)}+
\ldots +\xi_{\rho(n)})|^2 \mu(d\xi_1) \ldots \mu(d\xi_n).
\end{eqnarray*}

Proceeding as in the evaluation of $\psi^{*(n)}({\bf t}, {\bf t})$,
we obtain, %using Lemma \ref{hu-lemma},
\begin{eqnarray}
\nonumber \psi_{h}^{(n)}({\bf t},{\bf t}) &=& \int_{\bR^d}d\eta_1
\frac{\sin^2(u_1 |\eta_1|)}{|\eta_1|^2} |\eta_1|^{-\alpha}
\int_{\bR^d} d\eta_2
\frac{\sin^2(u_2|\eta_2|)}{|\eta_2|^2}|\eta_2-\eta_1|^{-\alpha}
\ldots \\
\label{formula-psi1} & & \int_{\bR^d}d\eta_n
\frac{|\sin((u_n+h)|\eta_n|) -\sin(u_n|\eta_n|)|^{2}}{|\eta_n|^2}
|\eta_n-\eta_{n-1}|^{-\alpha}.
\end{eqnarray}

By the Dominated Convergence Theorem,
$$\psi_{h}^{(n)}({\bf t},{\bf t}) \to 0 \quad \mbox{as} \ h \to 0.$$
The application of this theorem is justified, since
$$\frac{|\sin((t+h)|\xi|)-\sin(t|\xi|)|}{|\xi|} \leq \left(
\frac{4}{1+|\xi|^2}\right)^{1/2},$$ for all $\xi \in \bR^d, t>0, h
\in [0,1]$ (see p.4 of Erratum of \cite{dalang99}). The fact that
$E_1(t,h) \to 0$ follows by applying the Dominated Convergence
Theorem in (\ref{formula-E1}).

We now treat $E_2(t,h)$. Let $A=[0,t+h]^n \verb2\2 [0,t]^n$. We have
\begin{equation}
\label{formula-E2} E_2(t,h) =\alpha_H^n \int_{[0,t]^{2n}}
\prod_{j=1}^{n}|t_j-s_j|^{2H-2}1_{A}({\bf t}) 1_{A}({\bf s})
\gamma_{h}^{(n)}({\bf t}, {\bf s})d{\bf t} d{\bf s},
\end{equation}
where
$$\gamma_{h}^{(n)}({\bf t},{\bf s}) = \int_{\bR^{nd}} \cF g_{\bf
t}^{(n)}(\cdot \ , t+h,x)(\bxi)\overline{\cF g_{\bf s}^{(n)}(\cdot \
, t+h,x)(\bxi)} \mu(d\xi_1) \ldots \mu(d\xi_n).$$

By the Cauchy-Schwartz inequality, $\gamma_{h}^{(n)}({\bf t},{\bf
s}) \leq \gamma_{h}^{(n)}({\bf t},{\bf t})^{1/2}
\gamma_{h}^{(n)}({\bf s},{\bf s})^{1/2}$. To evaluate
$\gamma_h^{(n)}({\bf t},{\bf t})$, we use again
(\ref{Fourier-transform-gt}):
\begin{eqnarray}
\nonumber \gamma_{h}^{(n)}({\bf t},{\bf t})
% \|g_{\bf t}^{(n)}(\cdot \ ; t+h,x)\|_{\cP(\bR^d)^{\otimes n}}^{2} \\
& = & \int_{\bR^{nd}} |\cF g_{\bf t}^{(n)}(\cdot \ ,
t+h,x)(\bxi)|^2 \mu(d\xi_1) \ldots \mu(d\xi_n) \\
\nonumber &=& \int_{\bR^d}d\eta_1 \frac{\sin^2(u_1
|\eta_1|)}{|\eta_1|^2} |\eta_1|^{-\alpha} \int_{\bR^d} d\eta_2
\frac{\sin^2(u_2|\eta_2|)}{|\eta_2|^2}|\eta_2-\eta_1|^{-\alpha}
\ldots \\
\label{formula-psi2} & & \int_{\bR^d}d\eta_n
\frac{\sin^2((u_n+h)|\eta_n|)}{|\eta_n|^2}
|\eta_n-\eta_{n-1}|^{-\alpha},
\end{eqnarray}
where $u_j=t_{\rho(j+1)}-t_{\rho(j)}$. By Lemma \ref{hu-lemma},
$$\gamma_{h}^{(n)}({\bf t},{\bf t}) \leq
C_{\alpha,d}^n [u_1 \ldots u_{n-1}(u_n+h)]^{\alpha-d+2},$$ and hence
$$\gamma_{h}^{(n)}({\bf t},{\bf s}) \leq C_{\alpha,d}^n
\left[(u_n+h) (v_n+h)\prod_{j=1}^{n-1}u_jv_j
\right]^{(\alpha-d+2)/2},$$ where
$v_j=s_{\sigma(j+1)}-s_{\sigma(j)}$ and
$s_{\sigma(1)}<\ldots<s_{\sigma(n)}<s_{\sigma(n+1)}=t$.

By inequality (16) in \cite{balan-tudor09},
% the fact that $\|\varphi \|_{\cH(0,t+h)^{\otimes n}}^2 \leq b_H^{2n}
%\|\varphi \|_{L^{1/H}((0,t+h)^n)}^{2}$,
$$E_2(t,h) \leq b_H^{2n} C_{\alpha,d}^n \left(\int_{[0,t+h]^n}
1_{A}({\bf t})
\left[\prod_{j=1}^{n-1}(t_{\rho(j+1)}-t_{\rho(j)})(t+h-t_{\rho(n)})
\right]^{\delta} d{\bf t}\right)^{2H},$$ where
$\delta=(\alpha-d+2)/(2H)$. Using the fact that
$$A=\bigcup_{\rho \in S_{n-1}} \{(t_1, \ldots,t_n); 0<t_{\rho(1)}<
\ldots<t_{\rho(n-1)}<t_n \ \mbox{and} \ t_n \in [t,t+h]\},$$ we
obtain that:
\begin{eqnarray*}
E_2(t,h) & \leq &  b_H^{2n} C_{\alpha,d}^n \left[(n-1)! \int_t^{t+h}
(t+h-t_{n})^{\delta} I_{n-1}(t_n,\delta) dt_n \right]^{2H},
\end{eqnarray*}
where
\begin{eqnarray*}
I_{n-1}(t_n,\delta)&:=& \int_{0<t_{1}< \ldots < t_{n-1}<t_n}
\prod_{j=1}^{n-1}(t_{j+1}-t_{j})^{\delta }dt_1 \ldots dt_{n-1}\\
&=&
\frac{\Gamma(1+\delta)^{n}}{\Gamma((n-1)(1+\delta)+1)}t_n^{(n-1)(1+\delta)}
,\end{eqnarray*} (see Lemma 3.5 of \cite{balan-tudor09}). We obtain:
\begin{eqnarray*} E_2(t,h) & \leq & b_H^{2n} C_{\alpha,d}^n
\Gamma(1+\delta)^{2Hn}
\left[\frac{(n-1)!}{\Gamma((n-1)(1+\delta)+1)}\int_t^{t+h}
(t+h-t_n)^{\delta} t_n^{(n-1)(1+\delta)} dt_n \right]^{2H}\\
& \leq & C_{\alpha,d,H}^n(t+1)^{2H(n-1)(1+\delta)}
\frac{1}{[(n-1)!]^{2H \delta}}\left(\int_t^{t+h}(t+h-t_n)^{\delta}dt_n\right)^{2H} \\
& = & C_{\alpha,d,H}^n(t+1)^{2H(n-1)(1+\delta)}
\frac{1}{[(n-1)!]^{2H \delta}} h^{2H(\delta+1)} \to 0, \quad
\mbox{as} \ h \to 0.
%\left(\int_0^{h}u^{\delta}du\right)^{2H}
\end{eqnarray*}

%(For the heat equation, we use the same argument with
%$\delta=-(d-\alpha)/4H$.)

b) %Letting $y=x+z$,
We have:
\begin{eqnarray*}
E|J_n(t,x)-J_n(t,y)|^2&=& E|I_n(f_n (\cdot ,t,x))-I_n(f_n (\cdot
,t,y))|^2 =\frac{1}{n!}E_3(t,x,y),
\end{eqnarray*}
where
\begin{eqnarray}
\label{def-E3} E_3(t,x,y)&:=&(n!)^2 \|\tilde f_n (\cdot ,t,x)-\tilde
f_n (\cdot ,t,y) \|_{\cH
\cP^{\otimes n}}^{2} \\
\label{formula-E3} &=&\alpha_{H}^n
\int_{[0,t]^n}\prod_{j=1}^{n}|t_j-s_j|^{2H-2} \psi_{x,y}^{(n)}({\bf
t},{\bf s})d{\bf t} d{\bf s}
\end{eqnarray}
and
$$\psi_{x,y}^{(n)}({\bf t},{\bf s})=\int_{\bR^{nd}} \cF(g_{\bf
t}^{(n)}(\cdot \,t,x)-g_{\bf t}^{(n)}(\cdot \ ,t,y))(\bxi)
\overline{\cF(g_{\bf s}^{(n)}(\cdot \, t,x)-g_{\bf s}^{(n)}(\cdot \
,t,y))(\bxi)} \mu(d\xi_1) \ldots \mu(d\xi_n).$$

By the Cauchy-Schwartz inequality, $\psi_{x,y}^{(n)}({\bf t},{\bf
s}) \leq \psi_{x,y}^{(n)}({\bf t},{\bf t})^{1/2}
\psi_{x,y}^{(n)}({\bf s},{\bf s})^{1/2}$. To evaluate
$\psi_{x,y}^{(n)}({\bf t},{\bf t})$, we use
(\ref{Fourier-transform-gt}):
\begin{eqnarray}
\nonumber \psi_{x,y}^{(n)}({\bf t},{\bf t})
%\|g_{\bf t}^{(n)}(\cdot \, t,x)-g_{\bf t}^{(n)}(\cdot \ ,t,y)\|_{\cP(\bR^d)^{\otimes n}} \\
&=& \int_{\bR^{nd}} |\cF (g_{\bf t}^{(n)}(\cdot \ ,t,x)-g_{\bf
t}^{(n)}(\cdot \ ,t,y))(\bxi)|^2 \mu(d\xi_1) \ldots \mu(d\xi_n) \\
\nonumber  &=& \int_{\bR^d}d\eta_1 \frac{\sin^2(u_1
|\eta_1|)}{|\eta_1|^2} |\eta_1|^{-\alpha} \int_{\bR^d} d\eta_2
\frac{\sin^2(u_2|\eta_2|)}{|\eta_2|^2}|\eta_2-\eta_1|^{-\alpha}
\ldots \\
\label{formula-psi3} & & \int_{\bR^d}d\eta_n
\frac{\sin^2(u_n|\eta_n|)}{|\eta_n|^2} |\eta_n-\eta_{n-1}|^{-\alpha}
|1-e^{-i \eta_n \cdot (y-x)}|^2.
\end{eqnarray}
By the Dominated Convergence Theorem,
%$\psi_{x,y}^{(n)}({\bf t},{\bf t}) \to 0$ as $|x-y|\to 0$, and
$E_3(t,x,y) \to 0$ as $|x-y| \to 0$. $\Box$

\section{H\"older Continuity}

In this section, we obtain some bounds for the $p$-th moments of the
solution, from which we infer that the solution has a
$\gamma$-H\"{o}lder continuous modification, with
$0<\gamma<\frac{\alpha-d+2}{2}$.

If $f$ is the Riesz kernel of order $\alpha>d-2$, then for any
$1>\beta>(d-\alpha)/2$, $\int_{\bR^d}
\frac{\mu(d\xi)}{(1+|\xi|^2)^{\beta}}<\infty$, which is equivalent
to
\begin{equation}
\label{Conus-Dalang-relation} \sup_{\eta \in \bR^d} \int_{\bR^d}
\frac{\mu(d\xi)}{(1+|\xi+\eta|^2)^{\beta}}<\infty.
\end{equation}
(see (7.26) of \cite{conus-dalang09}).

%Assume that $\alpha>d-2$.
%Then we can find some $\beta \in
%(\frac{d-\alpha}{2},1)$ for which (\ref{Conus-Dalang-relation}) holds. Hence,
By Proposition 7.4 of \cite{conus-dalang09}, the fundamental
solution $G$ of the wave equation satisfies hypothesis (H3)-(H5) of
\cite{conus-dalang09}, for any $0<\gamma_i \leq
1-\beta<\frac{\alpha-d+2}{2}$, $i=1,2,3$. This fact is used in the
proof of the next result.

%{\bf Check if (H3)-(H5) hold for the heat equation as well.}

\begin{theorem}
Let $f$ be the Riesz kernel of order $\alpha>d-2$ and $u$ be the
solution of (\ref{wave}). Then for any $p \geq 2, T>0$ and $K
\subset \bR^d$ compact,
\begin{eqnarray*}
E|u(t+h,x)-u(t,x)|^p  & \leq & C |h|^{p[\gamma_1 \wedge
(\gamma_2+H)]}, \quad \forall t \in [0,T],\forall x \in \bR^d,
\forall h \in \bR,t+h \in
[0,T],\quad  \\
E|u(t,x+z)-u(t,x)|^p & \leq & C |z|^{p\gamma_3}, \quad \forall t
\geq 0,\forall x \in K,\forall z \in \bR^d, x+z \in K
\end{eqnarray*}
for any $0<\gamma_i<\frac{\alpha-d+2}{2}$, where $C$ is a constant
which depends on $\alpha,d,H,p,T$.

In particular, $\{u(t,x); (t,x) \in [0,T] \times K\}$ has a
modification
%$\{\tilde u(t,x); t \in [0,T], x \in K\}$
which is a.s. jointly $\gamma$-H\"older continuous in time and
space, for any $\gamma \in (0,\frac{\alpha-d+2}{2})$.
\end{theorem}

\noindent {\bf Proof:} We first treat the time increments. By
Minkowski's inequality and (\ref{norms-p-q-equiv}),
\begin{eqnarray}
\nonumber \lefteqn{\|u(t+h,x)-u(t,x) \|_{p} =\|\sum_{n \geq
0}(J_n(t+h,x)-J_n(t,x))
\|_{p} } \\
\nonumber & & \leq \sum_{n \geq 0}\|J_n(t+h,x)-J_n(t,x) \|_{p} \leq
 \sum_{n \geq 0} (p-1)^{n/2}\|J_n(t+h,x)-J_n(t,x) \|_{2} \\
\label{estimate-u-time}& & = \sum_{n \geq 0}
(p-1)^{n/2}\left\{\frac{2}{n!}(E_1(t,h)+E_2(t,h)) \right\}^{1/2},
\end{eqnarray}
where $E_1(t,h)$ and $E_2(t,h)$ are given by (\ref{def-E1}),
respectively (\ref{def-E2}).

To estimate $E_1(t,h)$, we use (\ref{formula-E1}) and
(\ref{formula-psi1}). The inner integral in (\ref{formula-psi1}) is
$$\int_{\bR^d}|\cF G(u_n+h, \cdot \ )(\xi+\eta_{n-1})-\cF G(u_n,
\cdot \ )(\xi+\eta_{n-1})|^2 \mu(d\xi).$$ This integral is bounded
by $C h^{2\gamma_1}$ for some $0<\gamma_1<\frac{\alpha-d+2}{2}$, due
to (H3). The remaining $(n-1)$-fold integral in (\ref{formula-psi1})
is bounded above by $C_{\alpha,d}^{n-1}(u_1 \ldots
u_{n-1})^{\alpha-d+2}$, by Lemma \ref{hu-lemma}. Hence,
$$\psi_{h}^{(n)}({\bf t},{\bf s}) \leq C h^{2\gamma_1}C_{\alpha,d}^{n-1}
(u_1 \ldots u_{n-1})^{(\alpha-d+2)/2}(v_1 \ldots
v_{n-1})^{(\alpha-d+2)/2},$$ where $u_j=t_{\rho(j+1)}-t_{\rho(j)}$
and $v_j=s_{\sigma(j+1)}-s_{\sigma(j)}$. By inequality (16) of
\cite{balan-tudor09},  % and Lemma 3.5 of \cite{balan-tudor09}.),
\begin{eqnarray}
\nonumber E_1(t,h) & \leq & C h^{2\gamma_1} C_{\alpha,d}^{n-1}
b_{H}^{2n} \left( n! \int_0^t \int_{0<t_1< \ldots <t_{n-1}<t_n}
\prod_{j=1}^{n-1}(t_j-t_{j-1})^{\delta} dt_1 \ldots dt_{n-1} dt_n
\right)^{2H} \\
\nonumber &=& C h^{2\gamma_1}
C_{\alpha,d}^{n-1}b_{H}^{2n}\Gamma(1+\delta)^{2Hn} \left(
\frac{n!}{\Gamma((n-1)(1+\delta)+1)} \int_0^t
t_n^{(n-1)(1+\delta)}dt_n \right)^{2H}, \\
\label{estimate-E2} & \leq & C h^{2\gamma_1} C_{\alpha,d,H}^{n}
\frac{T^{n(2H+\alpha-d+2)}}{(n!)^{\alpha-d+2}},
\end{eqnarray}
where $\delta=(\alpha-d+2)/(2H)$.

To estimate $E_2(t,h)$, we use (\ref{formula-E2}) and
(\ref{formula-psi2}). The inner integral in (\ref{formula-psi2}) is
$$\int_{\bR^d}|\cF G (u_n+h, \cdot \ )(\xi+\eta_{n-1})|^2\mu(d\xi) \leq
C(u_n+h)^{2\gamma_2},$$ for some $0<\gamma_2<\frac{\alpha-d+2}{2}$,
by (H4). Using the same ideas as above, we get:
\begin{eqnarray}
\nonumber E_2(t,h) & \leq & C C_{\alpha,d,H}^n
\left(\frac{(n-1)!}{\Gamma((n-1)(1+\delta)+1)}\int_t^{t+h}
(t+h-t_n)t_n^{(n-1)(1+\delta)}dt_n \right)^{2H} \\
\label{estimate-E1} & \leq & C h^{2(\gamma_2+H)}C_{\alpha,d,H}^{n}
\frac{T^{n(2H+\alpha-d+2)}}{(n!)^{\alpha-d+2}}.
\end{eqnarray}

From (\ref{estimate-u-time}), (\ref{estimate-E2}) and
(\ref{estimate-E1}), we get:
\begin{eqnarray*}
\|u(t+h,x)-u(t,x)\|_{p} & \leq & C h^{\gamma_1 \wedge
(\gamma_2+H)}\sum_{n \geq 0}\frac{(p-1)^{n/2}C_{\alpha,d,H}^{n/2}
T^{n(2H+\alpha-d+2)/2}}{(n!)^{(\alpha-d+3)/2}}\\
&=:& h^{\gamma_1 \wedge (\gamma_2+H)}C(\alpha,d,H,p,T).
\end{eqnarray*}

We now treat the spatial increments. As above, we obtain:
%By Minkowski's inequality and (\ref{norms-p-q-equiv}),
\begin{eqnarray*}
\|u(t,x+z)-u(t,x) \|_{p} & \leq &
% =\|\sum_{n \geq 0}(J_n(t,x+z)-J_n(t,x))\|_{p} } \\
%& & \leq \sum_{n \geq 0}\|J_n(t,x+z)-J_n(t,x) \|_{p} \leq
% \sum_{n \geq 0} (p-1)^{n/2}\|J_n(t,x+z)-J_n(t,x) \|_{2} \\
\sum_{n \geq 0} (p-1)^{n/2}\left(\frac{1}{n!}E_3(t,x,x+z)
\right)^{1/2},
\end{eqnarray*}
where $E_3(t,x,y)$ is given by (\ref{def-E3}). To estimate
$E_3(x,x+z)$ we use (\ref{formula-E3}) and (\ref{formula-psi3}).
Using the fact that $\cF G(u, \cdot-z)(\xi)=e^{-i \xi \cdot z} \cF
G(u,\cdot)(\xi)$, we see that the
 inner integral in (\ref{formula-psi3}) is:
$$\int_{\bR^d} |\cF G(u_n, \cdot -z )(\xi+\eta_{n-1})-
\cF G(u_n, \cdot \ )(\xi+\eta_{n-1})|^2 \mu(d\xi) \leq C |z|^{2
\gamma_3},$$ for some $0<\gamma_3<\frac{\alpha-d+2}{2}$, by (H5).
The rest of the proof is the same as above.

The final statement follows by a version of Kolmogorov's criterion
for multi-parameter processes (see e.g. Problem 2.9 of
\cite{karatzas-shreve91}). $\Box$

\section{Malliavin differentiability of the solution}
\label{Malliavin-differentiability}

In this section, we show that  $u(t,x)$ is Malliavin differentiable
of any order. When $d \leq 2$, we show that the Malliavin derivative
of the solution satisfies a certain integral equation. These results
are valid for the heat equation in any dimension $d$.

Recall that if $F$ is a smooth random variable of the form
(\ref{form-F}), the iterated Malliavin derivative $D^k F$ is an $\cH
\cP^{\otimes k}$-valued random variable, defined by:
$$D^k F=\sum_{i_1, \ldots,i_k=1}^{n}\frac{\partial^k f}{\partial x_{i_1} \ldots x_{i_k}}(W(\varphi_1), \ldots, W(\varphi_n))\varphi_{i_1} \otimes \ldots \otimes \varphi_{i_k}.$$

The space $\bD^{k,p}$ is the completion of the space $\cS$ of smooth
random variables, with respect to the the norm $\|\cdot
\|_{\bD^{k,p}}$ defined by:
\begin{equation}
\label{def-norm-Dkp} \|F\|_{\bD^{k,p}}^p=E|F|^p
+\sum_{j=1}^{k}\|D^{j}F\|_{\cH \cP^{\otimes j}}^p.
\end{equation}

%If $F=\sum_{n \geq 0}J_n F \in L^2(\Omega)$ and $k \geq 2$,
 %It is known that $F \in \bD^{k,2}$ if and only if $\sum_{n \geq k}n^{k}E|J_n F|^2<\infty$, and in this %case, $D^k (J_n F)=J_{n-k}F$ for any $n \geq k$, and
%$$E\|D^k F\|_{\cH \cP^{\otimes k}}^2=\sum_{n \geq k}n(n-1) \ldots (n-k+1)E|J_n F|^2$$
%(see p. 28 of \cite{nualart06}).

Using a known criterion (see e.g. p.28 of \cite{nualart06}), it
follows that $u(t,x) \in \bD^{k,2}$, since by (\ref{moment2-Jn}) and
(\ref{bound-alpha-n-wave}), we have:
$$\sum_{n \geq 1}n^kE|J_n(t,x)|^2=\sum_{n \geq 1}n^k\frac{1}{n!}\tilde \alpha_n(t)\leq \sum_{n \geq 1}
\frac{[2^kC(t)]^n}{(n!)^{\alpha-d+3}}<\infty.$$

Next, we show that $u(t,x) \in \bD^{k,p}$ for all $k \geq 1,p >1$.

%the fact that $n \leq 2^n$, we obtain:

%$$E\|D^k u(t,x)\|_{\cH \cP^{\otimes k}}^2=\sum_{n \geq k}n(n-1)
%\ldots (n-k+1)E|J_n(t,x)|^2$$

%\begin{equation}
%\label{sum-u-in-Dk} \sum_{n \geq 1}n^kE|J_n(t,x)|^2<\infty
%\end{equation}
%In the case of the wave equation,

%In the case of the heat
%equation, by (\ref{bound-alpha-n-heat}), we have:
%$$\sum_{n \geq 1}n^kE|J_n(t,x)|^2=\sum_{n \geq 1}n^k \frac{1}{n!}\alpha_n(t)\leq  \sum_{n \geq 1}
%\frac{[2^k C(t)]^n}{(n!)^{1-(d-\alpha)/2}}<\infty.$$

Let $F \in L^p(\Omega)$ be such that $D^k F$ exists. By Meyer's
inequalities (Theorem 1.5.1 of \cite{nualart06}), for any $p>1$,
\begin{equation}
\label{Meyer} E\|D^k F \|_{\cH \cP^{\otimes k}}^{p}<\infty \
\mbox{if and only if} \ E|C^k F|^p<\infty,
\end{equation}
%$$c_{k,p} E\|D^k F \|_{\cH^{\otimes k}}^{p} \leq E|C^k F|^p \leq
%C_{k,p} \{ E\|D^k F \|_{\cH^{\otimes k}}^{p}+ E|F|^p\},$$
where $C^k: {\rm Dom} \ C^k \subset L^2(\Omega) \to L^2(\Omega)$ is
the
operator defined by %$$CF=\sum_{n \geq 1}-\sqrt{n}J_n F.$$ Note that
$$C^k F =\sum_{n \geq 1}(-\sqrt{n})^k J_n F,$$
and ${\rm Dom} \ C^k=\{F \in L^2(\Omega); \sum_{n \geq 1}n^k E|J_n
F|^2<\infty\}=\bD^{k,2}$ for any $k \geq 1$.

By Minkowski's inequality and (\ref{norms-p-q-equiv}), we have:
\begin{equation}
\label{norm-C} \|C^k F\|_{p} \leq \sum_{n \geq 1} n^{k/2} \|J_n
F\|_{p} \leq \sum_{n \geq 1} n^{k/2} (p-1)^{n/2} \|J_n F\|_{2}.
\end{equation}

Combining (\ref{def-norm-Dkp}), (\ref{Meyer}) and (\ref{norm-C}), we
infer that $\|F\|_{\bD^{k,p}}<\infty$ (i.e. $F \in \bD^{k,p}$), if
\begin{equation}
\label{suf-cond-Dkp} \sum_{n \geq 1} n^{k/2} (p-1)^{n/2} \|J_n F
\|_{2}<\infty.
\end{equation}

Applying this in our case, we obtain the following result.

\begin{proposition}
Let $f$ be the Riesz kernel of order $\alpha>d-2$ and $u$ be the
solution of (\ref{wave}). Then $u(t,x) \in \bD^{k,p}$ for all $k
\geq 1$ and $p >1$.
\end{proposition}

\noindent {\bf Proof:} We verify (\ref{suf-cond-Dkp}). %In the case of the wave equation,
By (\ref{moment2-Jn}) and (\ref{bound-alpha-n-wave}), we have:
\begin{eqnarray*}
\sum_{n \geq 1}n^{k/2}(p-1)^{n/2}\|J_n(t,x)\|_2 &=& \sum_{n \geq
1}n^{k/2}(p-1)^{n/2} \left(\frac{1}{n!}\tilde \alpha_n(t)
\right)^{1/2}
\\ & \leq & \sum_{n \geq 1}
\frac{[2^k(p-1)C(t)]^{n/2}}{(n!)^{(\alpha-d+3)/2}}<\infty.
\end{eqnarray*}
%where we used the fact that $n \leq 2^n$ for all $n \geq 1$.
$\Box$

In the final part of this section, we show that the Malliavin
derivative $Du(t,x)$ satisfies a certain integral equation. For
this, we assume that $G(t,x)$ is a function in $x$ (i.e. $d \leq
2$).
%For this, we will assume that $d \leq 2$, and hence $G(t,x)$ is a
%function in $x$.

Recall that $u$ satisfies the integral equation
(\ref{u-integral-eq}). Intuitively, using the commutativity between
the operators $D$ and $\delta$, the derivative $Du$ should satisfy:
\begin{equation}
\label{anticipated-D-delta} Du(t,x)=G(t-\cdot,x-*)u+\int_0^t
\int_{\bR^d}G(t-s,x-y)Du(s,y)W(\delta s, \delta y),
\end{equation}
where $\cdot$ denotes the missing $r$ variable and $*$ denotes the
missing $z$ variable.

The integrand of the stochastic integral above is an $\cH \cP
\otimes \cH \cP$-valued random variable, and the integral needs to
be defined as an $\cH \cP$-valued random variable.

For this reason, we introduce a Hilbert-space-valued Skorohod
integral.

If $\cA$ is an arbitrary Hilbert space, we let $\cS(\cA)$ be the class of smooth $\cA$-valued random variables $F=\sum_{j=1}^{m}F_j v_j$, with $F_j \in \cS,v_j \in \cA, m \geq 1$. The Malliavin derivative of such $F$ is defined as $DF=\sum_{j=1}^{m}DF \otimes v_j$. % \in L^2(\Omega; \cH \cP \otimes V)$.
We denote by $\bD^{1,2}(\cA)$ the completion of $\cS(\cA)$ with
respect to the norm $\|\cdot \|_{\bD^{1,2}(\cA)}$, where
$$\|F\|_{\bD^{1,2}(\cA)}^{2}:=E\|F\|_{\cA}^2+E\|DF\|_{\cH \cP \otimes \cA}^2.$$

Similarly to the case $\cA=\bR$ (considered in Subsection \ref{Malliavin}),
we let $\delta^*$ be the adjoint of the operator $D$. %$D: \bD^{1,2}(\cA) \to L^2(\Omega;\cH \otimes\cA)$.
The domain of $\delta^*$, denoted by ${\rm Dom} \ \delta^*$, is the
set of $U \in L^2(\Omega;\cH \cP \otimes \cA)$ such that:
$$|E\langle DF,U \rangle_{\cH \cP \otimes \cA}| \leq c(E\|F\|_{\cA}^2)^{1/2}, \quad \forall F \in \bD^{1,2}(\cA),$$
where $c$ is a constant depending on $U$. If $U \in {\rm Dom} \
\delta^*$, then $\delta^*(U)$ is the element of $L^2(\Omega;\cA)$
characterized by the following duality relation:
\begin{equation}
\label{duality-delta*} E \langle F, \delta^*(U)
\rangle_{\cA}=E\langle DF,U \rangle_{\cH \cP \otimes \cA}, \quad
\forall F \in \bD^{1,2}(\cA).
\end{equation}

If $U\in {\rm Dom} \ \delta^*$, we use the notation
$$\delta^*(U)=\int_0^{\infty}\int_{\bR^d}U(t,x)W(\delta^* t,
\delta^* x),$$ even if $U(t,x)$ is not a function in $(t,x)$, and we
say that $\delta^*(U)$ is the $\cA$-valued Skorohod integral of $U$
with respect to $W$.

%If $D^{\varphi}F= \langle DF,\varphi \rangle_{\cH \cP}$ if $F \in
%\bD^{1,2}(\cH \cP)$ and $\varphi \in \cH \cP$.
Similar to the case $\cA=\bR$, we have the following result.

\begin{proposition}
\label{D-included-dom-delta} Let $U,V \in \bD^{1,2}(\cH \cP \otimes
\cA)$. Then
$$E \langle \delta^*(U), \delta^*(V)\rangle_{\cA}=
E \langle U,V\rangle_{\cH \cP \otimes
\cA}+E\left(\sum_{i,j,k=1}^{\infty} D^{e_i}\langle U,e_j \otimes a_k
\rangle_{\cH \cP \otimes \cA} D^{e_j}\langle V,e_i \otimes a_k
\rangle_{\cH \cP \otimes \cA}\right),$$ where $(e_i)_{i \geq 1}$,
$(a_k)_{k \geq 1}$ are complete orthonormal systems in $\cH \cP$,
respectively $\cA$. Consequently, if $U \in \bD^{1,2}(\cH \cP
\otimes \cA)$, then $U \in {\rm Dom} \ \delta^*$ and
\begin{equation}
\label{ineq-delta*}E\|\delta^*(U)\|_{\cA}^2 \leq E\|U \|_{\cH \cP
\otimes \cA}^2+E\|DU \|_{\cH \cP \otimes \cH \cP \otimes
\cA}^{2}=\|U\|_{\bD^{1,2}(\cH \cP \otimes \cA)}^2.
\end{equation}
\end{proposition}

\noindent {\bf Proof:} The proof is similar to Proposition 1.3.1 of
\cite{nualart06}. For this, one needs to revisit the basic rules of
the Malliavin calculus. %in the context of the Hilbert-space-valued Skorohod integral.
We omit the details, but we list below these rules, which apply to
any isonormal Gaussian process $\{W(h)\}_{h \in \cH}$:

1)  For any $F \in \cS(\cA), h \in \cH,v \in \cA$,
$$E\langle DF, h
\otimes v \rangle_{\cH \otimes \cA}=E (\langle F, v
\rangle_{\cA}W(h)).$$

2) For any $F \in \cS(\cA)$, $G \in \cS$, $h \in \cH$, $v \in \cA$,
$$E(G \langle DF,h \otimes v \rangle_{\cH \otimes \cA})=-E(\langle
F,v \rangle_{\cA} \langle DG,h \rangle_{\cH})+E( \langle F,v
\rangle_{\cA}GW(h)).$$

3) If $U=\sum_{j=1}^{m}F_j (h_j \otimes v_j)\in \cS(\cH \otimes
\cA)$ for some $F_j \in \cS$, $h_j \in \cH$, $v_j \in \cA$, then $U
\in {\rm Dom} \ \delta^*$ and
$$\delta^* (U)=\sum_{j=1}^{m}F_jW(h_j)-\sum_{j=1}^{m} \langle DF_j,
h_j\rangle_{\cH}v_j.$$

4) For any $U=\sum_{j=1}^{m}F_j (h_j \otimes v_j) \in \cS(\cH
\otimes \cA)$ and $h \in \cH$, $v \in \cA$,
$$D^{h \otimes v}(\delta^*(U))=\langle U,h \otimes v \rangle_{\cH \otimes \cA}
+\langle \delta^*(D^h u),v \rangle_{\cA},$$ where $D^{h \otimes
v}(\delta^*(U))=\langle D(\delta^*(U)),h \otimes v \rangle_{\cH
\otimes \cA}$ and $D^h U=\sum_{j=1}^{m}(D^h F_j) (h_j \otimes v_j)$.
Here $D^h F=\langle DF,h \rangle_{\cH}$. $\Box$

\vspace{3mm}

In what follows, we let $\delta^*$ be the operator corresponding to
the case $\cA=\cH \cP$.

We begin with some preliminary results.

\begin{lemma}
\label{bound-on-D-un}For any $t>0,x \in \bR^d$, $Du_n(t,x) \to Du(t,x)$ in %the weak topology of
$L^2(\Omega;\cH \cP)$, and
\begin{equation}
\label{D-bounded-in-norm} C_T^{(1)}:= \sup_{(t,x) \in [0,T] \times
\bR^d} E\|Du(t,x)\|_{\cH \cP}^2<\infty, \quad \mbox{for all} \ T>0.
\end{equation}
\end{lemma}

\noindent {\bf Proof:} Using Lemma 1.2.3 of \cite{nualart06}, %with $F_n=u_n(t,x)$ and $F=u(t,x)$,
it suffices to prove that:
$$\sup_{n \geq 1}\sup_{(t,x) \in [0,T] \times \bR^d} E\|Du_n(t,x)\|_{\cH \cP}^2<\infty.$$

%By abuse of notation, we write $D_{r,z}u_n(t,x)$ even if this is not
%a function in $(r,z)$.
By Proposition \ref{Malliavin-rule-1},
$D_{r,z}u_n(t,x)%=\sum_{k=1}^{n}D_{r,z}J_k(t,x)=\sum_{k=1}^{n}D_{r,z}(I_k(f_k(\cdot \,t,x)))
=\sum_{k=1}^{n}k I_{k-1}(f_k(\cdot \ , r,z,t,x))$.
% where $\bullet$ denotes the missing $k$-th variable of $f_k(\cdot,t,x)$.
Using the orthogonality of the Wiener chaos spaces,
(\ref{isometry-In}) and (\ref{S-finite-2-strong}), we get: for any
$t \in [0,T]$
%\begin{eqnarray*}
%E\|Du_n(t,x)\|_{\cH \cP}^2 &=& \sum_{k=1}^{n}k^2 E \|
%I_{k-1}(f_k(\cdot \ ,\bullet, t,x)) \|_{\cH \cP}^2 \\
%& \leq & \sum_{k=1}^{n}k^2 (k-1)!  \|f_k(\cdot \ ,
%t,x)\|_{\cH \cP^{\otimes k}}^2 \\
%\langle \langle \tilde f_k(\cdot \ , \bullet,t,x),\tilde f_k(\cdot \ \bullet,t,x) \rangle_{\cH %\cP^{\otimes(k-1)}} \rangle_{\cH \cP} \\
%&=&  \sum_{k=1}^{n}k k! \ \|\tilde f_k(\cdot \ , t,x)\|_{\cH
%\cP^{\otimes k}}^2 \\
%&=& \sum_{k=1}^{n}k k! \frac{1}{(k!)^2}\alpha_k(t)=\sum_{k=1}^{n}
%\frac{1}{(k-1)!}\alpha_k(t) \leq C_T<\infty,
%\end{eqnarray*}
%for any $t \in [0,T],x \in \bR^d,n \geq 1$. $\Box$
\begin{eqnarray*}
E\|Du_n(t,x)\|_{\cH \cP}^2
 &=& \alpha_{H} \int_{(0,t)^2}\int_{\bR^{2d}}|r-r'|^{2H-2} f(z-z')
\left(\sum_{k=1}^{n}k I_{k-1}(f_k(\cdot \ , r,z,t,x)) \right) \\
& &  \left(\sum_{l=1}^{n}l I_{l-1}(f_l(\cdot \ , r',z',t,x))
 \right)dz dz' dr dr' \\
&=& \alpha_{H} \sum_{k=1}^{n}k^2 (k-1)!  \int_{(0,t)^2}\int_{\bR^{2d}}|r-r'|^{2H-2} f(z-z') \\
& & \langle \bar f_k(\cdot \ , r,z,t,x),\bar f_k(\cdot \ ,r',z',t,x)
 \rangle_{\cH \cP^{\otimes(k-1)}}dz dz' dr dr' \\
 &=&  \sum_{k=1}^{n}k k! \ \|\bar f_k(\cdot \ , t,x)\|_{\cH
 \cP^{\otimes k}}^2
 \leq \sum_{k=1}^{n}k k! \ \|f_k(\cdot \ , t,x)\|_{\cH
 \cP^{\otimes k}}^2 \\
 &=& \sum_{k=1}^{n}k k! \frac{1}{(k!)^2}\alpha_k(t)=\sum_{k=1}^{n}
 \frac{1}{(k-1)!}\alpha_k(t) \leq C_T<\infty,
 \end{eqnarray*}
where $\bar f_{k}(\cdot,r,z,t,x)$ denotes the symmetrization of
$f_k(\cdot,r,z,t,x)$ with respect to the first $k-1$ variables. For
the first inequality above, we used the fact that the $\|\tilde f
\|_{\cH \cP^{\otimes n}} \leq \|f\|_{\cH \cP^{\otimes n}}$ for any
$f \in \cH \cP^{\otimes n}$. For the last inequality, we used Remark
\ref{remark-alpha-non-sym}. $\Box$

\begin{remark}
{\rm Using Proposition \ref{Malliavin-rule-1} iteratively, we obtain
that:
$$D_{(\tau,w),(r,z)}^2u(t,x)=\sum_{n \geq
2}n(n-1)I_{n-2}(f_n(\cdot,\tau,w,r,z,t,x)),$$ for any $(\tau,w) \in
(0,t) \times \bR^d$ and $(r,z) \in (0,t) \times \bR^d$. Hence,
similarly to (\ref{D-bounded-in-norm}), one can show that:
\begin{equation}
\label{D2-bounded-in-norm} C_T^{(2)}:= \sup_{(t,x) \in [0,T] \times
\bR^d} E\|D^2u(t,x)\|_{\cH \cP \otimes \cH \cP}^2<\infty, \quad
\mbox{for all} \ T>0.
\end{equation}
}
\end{remark}

\begin{lemma}
\label{U-belongs-dom} For any $t>0, x\in \bR^d$, the process
$U^{(t,x)}$  defined by
\begin{equation}
\label{def-U-tx} U^{(t,x)}=\{U^{(t,x)}(s,y):=G(t-s,x-y)Du(s,y);  s
\geq 0, y \in \bR^d\}
\end{equation}
belongs to ${\rm Dom} \ \delta^*$.
\end{lemma}

\noindent {\bf Proof:} By Proposition \ref{D-included-dom-delta}, it
suffices to show that $U^{(t,x)} \in \bD^{1,2}(\cH \cP \otimes \cH
\cP)$, i.e. $\|U^{(t,x)} \|_{\bD^{1,2}(\cH \cP \otimes \cH
\cP)}<\infty$. Note that
$$\|U^{(t,x)} \|_{\bD^{1,2}(\cH \cP \otimes \cH \cP)}^2=
E\|U^{(t,x)}\|_{\cH \cP \otimes \cH \cP}^2+E\|DU^{(t,x)}\|_{\cH \cP
\otimes \cH \cP \otimes \cH \cP}^2.$$

%By Proposition \ref{Malliavin-rule-1}, $Du (s,y)=\sum_{n \geq
%1}nI_{n-1}(f_n(\cdot,\bullet,s,y))$ and $D^2u(s,y)=\sum_{n \geq
%2}n(n-1)I_{n-2}(f_n(\cdot,\bullet,\bullet,s,y))$. Using the same
%argument as in the proof of Theorem \ref{prop-existence}, one can
%show that
%\begin{eqnarray*}
%U^{(t,x)}&=&\sum_{n \geq
%1}nI_{n-1}(f_{n+1}(\cdot,\bullet,\bullet,t,x)) \\
%DU^{(t,x)}&=&\sum_{n \geq 1}n(n-1)I_{n-2}(f_{n+1}(\cdot,
%\bullet,\bullet,\bullet,t,x)).
%\end{eqnarray*}
% Hence,
%\begin{eqnarray*}
%E\|U^{(t,x)}\|_{\cH \cP \otimes \cH \cP}^{2}&=& \sum_{n \geq 1}n^2 E
%\| I_{n-1}(f_{n+1}(\cdot,\bullet,\bullet,t,x))\|_{\cH \cP \otimes
%\cH \cP}^2 \\
%& \leq & \sum_{n \geq 1}n^2 (n-1)!\|f_{n+1}(\cdot,t,x)\|_{\cH
%\cP^{\otimes (n+1)}}^2 \\
%&=& \sum_{n \geq 1}n n! \frac{1}{[(n+1)!]^2}\alpha_{n+1}(t)<\infty,
%\end{eqnarray*}
%and
%\begin{eqnarray*}
%E\|DU^{(t,x)}\|_{\cH \cP^{\otimes 3}}^{2}&=& \sum_{n \geq
%2}n^2(n-1)^2 E \|
%I_{n-2}(f_{n+1}(\cdot,\bullet,\bullet,\bullet,t,x))\|_{\cH \cP^{\otimes 3}}^2 \\
%& \leq & \sum_{n \geq 2}n^2(n-1)^2 (n-2)!\|f_{n+1}(\cdot,t,x)\|_{\cH
%\cP^{\otimes (n+1)}}^2 \\
%&=& \sum_{n \geq 2}n(n-1) n!
%\frac{1}{[(n+1)!]^2}\alpha_{n+1}(t)<\infty.
%\end{eqnarray*}
By the Cauchy-Schwartz inequality, (\ref{D-bounded-in-norm}) and
(\ref{D2-bounded-in-norm}),
\begin{eqnarray*}
E\|U\|_{\cH \cP^{\otimes 2}}^2 &=& \alpha_H \int_0^t \int_0^t
G(t-s,x-y)G(t-s',x-y)
E\langle Du(s,y),Du(s',y') \rangle_{\cH \cP} \\
& & |s-s'|^{2H-2}f(y-y') dy dy' ds ds' \\
& \leq & C_t^{(1)} \|G(t-\cdot,x-*)\|_{\cH \cP}^2<\infty,  \\
E\|DU\|_{\cH \cP^{\otimes 3}}^2 &=& \alpha_H \int_0^t \int_0^t
G(t-s,x-y)G(t-s',x-y)
E\langle D^2u(s,y),D^2u(s',y') \rangle_{\cH \cP^{\otimes 2}} \\
& & |s-s'|^{2H-2}f(y-y') dy dy' ds ds' \\
& \leq & C_t^{(2)} \|G(t-\cdot,x-*)\|_{\cH \cP}^2<\infty.
\end{eqnarray*}
$\Box$

Using the same argument as above, %in Lemma \ref{U-belongs-dom},
one can show that the process $U_{n}^{(t,x)}$ defined by
$U_{n}^{(t,x)}(s,y):=G(t-s,x-y)Du_{n}(s,y)$, belongs to ${\rm Dom} \
\delta^*$. %, for any $n \geq 0$.

The next result shows that the sequence $\{Du_{n}(t,x)\}_{n \geq 0}$
satisfies a recurrence relation. %, in the spirit of (\ref{anticipated-D-delta}).

\begin{proposition}
\label{equation-Du-n} For any $t>0,x \in  \bR^d$ and $n \geq 1$,
$$Du_{n}(t,x)=G(t-\cdot,x-*)u_{n-1}+\int_0^t \int_{\bR^d}U_{n-1}^{(t,x)}(s,y)W(\delta^* s, \delta^* y)$$
in $L^2(\Omega; \cH \cP)$.
\end{proposition}

\noindent {\bf Proof:}
%Let $\cP(\bR^d)$ be the Hilbert space,
%defined as the closure of the space $\cE(\bR^d)$ of linear
%combinations of elementary functions on $\bR^d$, with respect to the
%inner product $\langle 1_{A},1_{B} \rangle_{\cP(\bR^d)}=\int_A
%\int_B f(x-y)dx dy$. Note that $\cP(\bR^d)$ may contain
%distributions in $\cS'(\bR^d)$.
%In the argument below, we denote by $*$ the missing $z$-variable,
%of the elements in $\cP(\bR^d)$,
%and we use the formal notation (\ref{formal-def-fn}) for the kernels
%$f_n(\cdot,t,x)$.
{\em Step 1.} By induction on $n$, we show that for any $r \in
(0,t)$, $z \in \bR^d$,
$D_{r,z}u_{n}(t,x)=\sum_{k=1}^{n}A_k$, %in $\cP(\bR^d)$,
where $A_k=\sum_{i=1}^{k}A_k^{(i)}$,
\begin{eqnarray*}
A_{k}^{(1)}&=& \int_{r<s_1}G(t-s_{k-1},x-y_{k-1}) \ldots
G(s_1-r,y_1-z)W(ds_1,dy_1)
\ldots W(ds_{k-1}, dy_{k-1}) \\
A_{k}^{(i)}&=& \int_{s_{i-1}<r<s_i}G(t-s_{k-1},x-y_{k-1}) \ldots
G(s_i-r,y_i-z)G(r-s_{i-1},z-y_{i-1})  \\
& & \ldots G(s_2-s_1,y_2-y_1)W(ds_1,dy_1)
\ldots W(ds_{k-1}, dy_{k-1}), \quad i=2, \ldots,k-1\\
A_{k}^{(k)}&=& \int_{s_{k-1}<r}G(t-r,x-z)G(r-s_{k-1},z-y_{k-1})
\ldots G(s_2-s_1,y_2-y_1)\\
& & W(ds_1,dy_1) \ldots W(ds_{k-1}, dy_{k-1}),
\end{eqnarray*}
and the integrals above are taken over the set
$\{0<s_1<\ldots<s_{k-1}<t\} \times \bR^{nd}$.

By definition, $u_{n}(t,x)=u_{n-1}(t,x)+I_n(\tilde f_n(\cdot,t,x))$.
Using Proposition \ref{Malliavin-rule-1} and the induction
hypothesis, we obtain that:
$$D_{r,z}u_{n}(t,x)=%D_{r,z}u_{n-1}(t,x)+nI_{n-1}(f_n(\cdot,r,z,t,x))=
\sum_{k=1}^{n-1}A_k+nI_{n-1}(\tilde f_n(\cdot,r,z,t,x)).$$ Note that
$nI_{n-1}(\tilde f_n(\cdot,r,z,t,x))=A_n$, the $n$ integrals
$A_1^{(n)}, \ldots, A_n^{(n)}$ corresponding to the $n$ possible
locations of $r$, compared with the variables $s_1<\ldots<s_{n-1}$.

{\em Step 2.} We prove that for every $r \in (0,t)$, $z \in \bR^d$
$$D_{r,z}u_n(t,x)=G(t-r,x-z)u_{n-1}(r,z)+\int_0^t \int_{\bR^d}
G(t-s,x-y)D_{r,z}u_{n-1}(s,y)W(\delta s,\delta y).$$
% in $\cP(\bR^d)$.

We use the expression of $D_{r,z}u_n(t,x)$ obtained in Step 1. Note
that the terms $A_1^{(1)},A_2^{(2)}, \ldots,A_n^{(n)}$ have the
common factor $G(t-r,x-z)$. A quick calculation shows that the sum
of these terms is $G(t-r,x-z)u_{n-1}(r,z)$.
%\begin{eqnarray*}
%\lefteqn{G(t-r,x-y)\left[1+ \int_{s_1<r} G(r-s_1,z-y_1)W(ds_1,dy_1)+ \ldots \right. } \\
%& & \left.\int_{s_1< \ldots<s_{n-1}<r}G(r-s_{n-1},z-y_{n-1}) \ldots G(s_2-s_1,y_2-y_1)W(ds_1,dy_1) \ldots
%W(ds_{n-1},dy_{n-1}) \right] \\
%& & =G(t-r,x-y)u_{n-1}(r,z).
%\end{eqnarray*}

For the remaining terms, we change the names of the variables of
integration, so that $G(t-s,x-y)$ becomes a common factor. More
precisely, we call $(s,y)$ the variable $(s_{k-1},y_{k-1})$ in
$A_k^{(i)}$, for any $i=1,\ldots, k-1$ and $k=2, \ldots,n$. The sum
of these terms turns out to be
$$\int_0^t \int_{\bR^d} G(t-s,x-y)D_{r,z}u_{n-1}(s,y)W(ds,dy),$$
using the expression of $D_{r,z}u_{n-1}(s,y)$ obtained in Step 1.
Finally, by (\ref{replace-d-by-delta}), %Proposition \ref{Malliavin-rule-2},
we can replace the integral $W(ds,dy)$ by an integral $W(\delta
s,\delta y)$.

{\em Step 3.} We show that the process $\delta^*(U_{n-1}^{(t,x)})$
coincides (in $L^2(\Omega,\cH \cP)$) with the process
$V_{n-1}^{(t,x)}$, defined by
\begin{eqnarray*}V_{n-1}^{(t,x)}(r,z)&:=& \int_0^t \int_{\bR^d}G(t-s,x-y)
D_{r,z}u_{n-1}(s,y)W(\delta s,\delta y) \\
&=& \delta(G(t-\cdot,x-*) D_{r,z}u_{n-1}).
\end{eqnarray*}
%where $\bullet$ denotes the missing $(r,z)$-variable.

By the duality relation (\ref{duality-delta*}), it suffices to prove
that
\begin{equation}
\label{duality-U-V} E\langle F,V_{n-1}^{(t,x)} \rangle_{\cH \cP}=E
\langle DF,U_{n-1}^{(t,x)}\rangle_{\cH \cP \otimes \cH \cP}, \quad
\forall F \in \bD^{1,2}(\cH \cP).
\end{equation}

Without loss of generality, we may assume that $F$ is smooth, i.e.
$F=F_0 \varphi$ with $F_0 \in \cS$ and $\varphi \in \cH \cP$. Then
$DF=DF_0 \otimes \varphi$ and
%\begin{eqnarray*}
%E \langle DF,U_{n-1}^{(t,x)}\rangle_{\cH \cP \otimes \cH \cP}
 %E \langle  \langle DF_{0}, G(t-*,x-\cdot) D_{\bullet}u_n \rangle_{\cH \cP} , \ \varphi(\bullet)
 % \rangle_{\cH \cP}  \\
%&=& \langle  E \langle DF_{0}, G(t-*,x-\cdot) D_{\bullet}u_n \rangle_{\cH \cP} , \ \varphi(\bullet) \rangle_{\cH \cP} \\
%&=& \langle E(F_{0}V_{n-1}^{(t,x)}(\bullet))  , \ \varphi(\bullet) \rangle_{\cH \cP} \\
%&=& E \langle F_{0}V_{n-1}^{(t,x)}(\bullet)  , \ \varphi(\bullet) \rangle_{\cH \cP} \\
%&=& E\langle F, V^{(t,x)} \rangle_{\cH \cP}.
%\end{eqnarray*}
%which in formal notation can be written as:
\begin{eqnarray*}
\lefteqn{E \langle DF,U_{n-1}^{(t,x)}\rangle_{\cH \cP \otimes \cH \cP} = \alpha_H \int_{(\bR_{+} \times \bR^d)^2} E \langle DF_{0}, G(t-*,x-\cdot) D_{r,z}u_n \rangle_{\cH \cP} } \\
& & \quad \varphi(r',z')|r-r'|^{2H-2}f(z-z')dz dz' dr dr' \\
& & =\alpha_H \int_{(\bR_{+} \times \bR^d)^2} E(F_{0}V_{n-1}^{(t,x)}(r,z))  \varphi(r',z')|r-r'|^{2H-2}f(z-z')dz dz' dr dr' \\
& & = E\langle F, V^{(t,x)} \rangle_{\cH \cP}.
\end{eqnarray*}
Note that for the second-last equality above, we used the duality
relation (\ref{duality}) (for the operator $\delta$), whereas for
the first and last equality we used Fubini's theorem. This shows
(\ref{duality-U-V}), and concludes the proof. $\Box$

\vspace{3mm}

The following result gives the precise meaning of relation (\ref{anticipated-D-delta}).% mentioned earlier.

\begin{theorem}
Let $f$ be the Riesz kernel of order $\alpha>d-2$ and $u$ be a
solution of (\ref{wave}). For any $t>0,x\in \bR^d$, let $U^{(t,x)}$
be defined by (\ref{def-U-tx}). Then,
$$Du(t,x)=G(t-\cdot,x-*)u+\int_0^t \int_{\bR^d}U^{(t,x)}(s,y)W(\delta^* s, \delta^* y)$$
in $L^2(\Omega; \cH \cP)$.
\end{theorem}

\noindent {\bf Proof:} By the duality relation
(\ref{duality-delta*}), it suffices to prove that,
\begin{equation}
\label{identity-Malliavin}E \langle
Du(t,x)-G(t-\cdot,x-*)u,F\rangle_{\cH \cP}= E\langle DF, U^{(t,x)}
\rangle_{\cH \cP \otimes \cH \cP},
\end{equation}  for any $F \in \bD^{1,2}(\cH \cP)$.
By Proposition \ref{equation-Du-n}, for any $F \in \bD^{1,2}(\cH
\cP)$,
\begin{equation}
\label{identity-Malliavin-n} E \langle
Du_n(t,x)-G(t-\cdot,x-*)u_{n-1},F\rangle_{\cH \cP}= E\langle DF,
U_{n-1}^{(t,x)} \rangle_{\cH \cP \otimes \cH \cP}.
\end{equation}

Relation (\ref{identity-Malliavin}) is obtained by taking $n \to
\infty$ in (\ref{identity-Malliavin-n}). We justify this below.

On the right-hand side of (\ref{identity-Malliavin-n}), we use the
duality relation (\ref{duality-delta*}), the Cauchy-Schwartz
inequality, and (\ref{ineq-delta*}):
\begin{eqnarray*}
E \langle DF, U_{n-1}^{(t,x)}-U^{(t,x)} \rangle_{\cH \cP \otimes \cH \cP} &=& E \langle \delta^*(U_{n-1}^{(t,x)}-U^{(t,x)}),F \rangle_{\cH \cP}  \\
& \leq & (E \|\delta^*(U_{n-1}^{(t,x)}-U^{(t,x)})|_{\cH \cP}^2)^{1/2} (E\|F \|_{\cH \cP}^2)^{1/2} \\
& \leq & \|U_{n-1}^{(t,x)}-U^{(t,x)}\|_{\bD^{1,2}(\cH \cP \otimes
\cH \cP)}^2 (E\|F \|_{\cH \cP}^2)^{1/2}.
\end{eqnarray*}
To show that $\|U_{n-1}^{(t,x)}-U^{(t,x)}\|_{\bD^{1,2}(\cH \cP
\otimes \cH \cP)}^2 \to 0$ as $n \to \infty$, we use the same
argument as in Lemma \ref{U-belongs-dom}.

For the first term on the left-hand side of
(\ref{identity-Malliavin-n}), by the Cauchy-Schwartz inequality and
Lemma \ref{bound-on-D-un},
$$E \langle Du_n(t,x)-Du(t,x),F\rangle_{\cH \cP} \leq
(E\|Du_n(t,x)-Du(t,x)\|_{\cH \cP}^2)^{1/2}(E\|F\|_{\cH \cP}^2)^{1/2}
\to 0,$$ as $n \to \infty$. For the second term on the left-hand
side of (\ref{identity-Malliavin-n}), we assume that $F$ is smooth,
i.e. $F=F_0 \varphi$ with $F_0 \in \cS$ and $\varphi \in \cH \cP$.
%Let $G_{\varepsilon}(t,\cdot)=\psi_{\varepsilon}*G(t,\cdot)$ be a
%regularization of $G(t,\cdot)$, i.e
%$\psi_{\varepsilon}(x)=\varepsilon^{-d}\psi(x/\varepsilon)$, $\psi
%\in C_0^{\infty}(\bR^d)$, $\psi \geq 0$, $\int_{\bR^d} \psi(x)dx=1$.
By the Cauchy-Schwartz inequality and the Dominated Convergence
theorem,
\begin{eqnarray*}
\lefteqn{E \langle G(t-\cdot,x-*)(u_{n-1}-u),F \rangle_{\cH \cP} =
\alpha_H \int_{(\bR_{+} \times \bR^d)^2} E[F_0(u_{n-1}(s,y)-u(s,y))]
} \\ & & G(t-s,x-y) \varphi(s',y') |s-s'|^{2H-2}f(y-y')dy dy ds ds'
\\
& & \leq \alpha_H E(F_0^2)^{1/2} \int_{(\bR_{+} \times \bR^d)^2}
(E|u_{n-1}(s,y)-u(s,y))|^2)^{1/2} G(t-s,x-y) \varphi(s',y')\\
& & \quad |s-s'|^{2H-2}f(y-y')dy dy ds ds', \quad \mbox{as $n \to
\infty$.}
\end{eqnarray*}
%uniformly in $\varepsilon>0$. Taking $\varepsilon \to 0$, we infer
%that $\lim_{n \to \infty}E \langle G(t-\cdot,x-*) \linebreak
%(u_{n-1}-u),F \rangle_{\cH \cP} = 0$.
This concludes the proof of
(\ref{identity-Malliavin}). $\Box$

\end{document}